\documentclass[12pt]{article}
\usepackage{epsfig}
\usepackage{amsmath}
\usepackage{amssymb}
\usepackage{eucal}

\newtheorem{thm}{Theorem}[section]

\newtheorem{prop}[thm]{Proposition} 
\newtheorem{rmk}[thm]{Remark}

\newtheorem{hyp}[thm]{Hypothesis}

\newcommand{\proof}{\noindent {\bf Proof} \hspace{0.2in}} 
\newcommand{\qed}{\hfill\mbox{\raggedright\rule{.07in}{.1in}}
  \vspace{1ex}} 

\newcommand{\Section}[1]{\section{#1} \setcounter{equation}{0}}

\title{Zero-Hopf bifurcation in the Van der Pol oscillator with delayed position and velocity feedback}

 \author{Jason Bramburger \and
Benoit Dionne \and
Victor G. LeBlanc\\Department of
  Mathematics and Statistics\\University of Ottawa\\Ottawa, 
ON K1N 6N5\\CANADA}
\date{\today}

\begin{document}

\maketitle

\begin{abstract}
In this paper, we consider the traditional Van der Pol Oscillator with a forcing dependent on a delay in feedback. The delay is taken to be a nonlinear function of both position and velocity which gives rise to many different types of bifurcations. In particular, we study the Zero-Hopf bifurcation that takes place at certain parameter values using methods of centre manifold reduction of DDEs and normal form theory. We present numerical simulations that have been accurately predicted by the phase portraits in the Zero-Hopf bifurcation to confirm our numerical results and provide a physical understanding of the oscillator with the delay in feedback. 
\end{abstract}

\pagebreak
\Section{Introduction}

The classical van der Pol oscillator
\begin{equation} \label{OrigVanderPol}
	\ddot{x}(t) + \varepsilon(x^2(t) - 1)\dot{x}(t) + x(t) = f(t),	
\end{equation}
is one of the most well-studied paradigms for nonlinear oscillators. Originally equation (\ref{OrigVanderPol}) was a model for an electrical circuit with a triode valve, and it has evolved into one of the most celebrated equations in the study of nonlinear dynamics.  Typically $\varepsilon >0$ is a parameter, and $f(t)$ is the external forcing term.

It is well-known that in the unforced system (when $f \equiv 0$) the trivial equilibrium point is unstable and the system contains a stable limit cycle oscillation for all values of $\varepsilon$. 
In the case of time-periodic forcing, more complicated dynamics can occur, and even chaotic behaviour is possible if $\varepsilon$ is large enough.

The van der Pol oscillator has become synonymous in many far reaching branches of the physical sciences with modern systems that exhibit limit cycle oscillations. Kaplan et al. provide a simple example of a biological application of the van der Pol oscillator to the dynamics of the heart \cite{Kaplan}, and the FitzHugh-Nagumo model is a modified form of the oscillator that has direct applications to neurones in the brain. On the larger scale, Cartwright et al. showed that the van der Pol oscillator can be used to model earthquake faults with viscous friction \cite{Cartwright}. All of these applications and many more make the van der Pol oscillator one of the most widely celebrated differential equations in mathematics and one that is just as relevant today as the day that van der Pol himself presented the equation.  

Despite all of these applications of the van der Pol oscillator, it wasn't until around the start of the twenty-first century, almost eighty years after van der Pol introduced his equation, that mathematicians began to investigate the effects of delayed feedback on the equation.
Delay-differential equations are used as models in many areas of science, engineering, economics and beyond \cite{BBL, HFEKGG,Kuang,LM,SieberKrauskopf,SC,SS,VTK,
WuWang,ZW}.
It is now well understood that retarded functional differential equations (RFDEs), a class which contains delay-differential equations, behave for the most part like ordinary differential equations on appropriate infinite-dimensional function spaces.  As such, many of the techniques and theoretical results of finite-dimensional dynamical systems have counterparts in the theory of RFDEs.  In particular, versions of the stable/unstable and center manifold theorems in neighborhoods of an equilibrium point exist for RFDEs \cite{HVL}.  Also, techniques for simplifying vector fields via center manifold and normal form reductions have been adapted to the study of bifurcations in RFDEs \cite{FM1,FM2}.
 
Atay \cite{Atay} studied the equation
\begin{equation}
	\ddot{x}(t) + \varepsilon(x^2(t) - 1)\dot{x}(t) + x(t) = \varepsilon kx(t - \tau),	
\end{equation} 
where $\tau > 0$ is a delay.  In particular, he showed that the delay can change the stability of the limit cycle which is present in the non-delayed version  (\ref{OrigVanderPol}).  A similar study was presented by Oliveira \cite{deOlivera}, in which delay terms were added to the nonlinear part of the unforced system, and found that periodic solutions still existed under certain conditions. 

Wei and Jiang \cite{WeiJiang,WeiJiang2} study the general case of a delayed position feedback forcing of the van de Pol equation
\begin{equation}
	\ddot{x}(t) + \varepsilon(x^2(t) - 1)\dot{x}(t) + x(t) = g(x(t-\tau)).	
	\label{eq1}
\end{equation}
They found that under certain conditions single, double and triple zero eigenvalues were possible as well as a purely imaginary pair with a zero eigenvalue. The cases of the single and double zero eigenvalues were presented in \cite{WeiJiang2} and \cite{JiangYuan}, respectively, with full details of the nature of their respective bifurcations. Later, Wu and Wang \cite{WuWang} undertook the bifurcation analysis of the purely imaginary pair and the zero eigenvalue to find that a Zero-Hopf bifurcation was occurring. They found that only two of the four possible bifurcation diagrams of the Zero-Hopf bifurcation were possible.  This is due to the nonlinear forcing term $g$ in (\ref{eq1}) that depends only on delayed position.  In principle, it could also depend on delayed velocity. 

In this paper, we consider the more general case of a feedback forcing which is dependent on both delayed position and delayed velocity
\begin{equation} \label{VanDerPol}
	\ddot{x}(t) + \varepsilon(x^2(t) - 1)\dot{x}(t) + x(t) = g(\dot{x}(t-\tau),x(t-\tau)),
\end{equation}
where g $\in C^3$, $g(0,0) = 0$, $g_{\dot{x}}(0,0) = a$ and $g_x(0,0) = b$, and work to understand its dynamics. This allows us to consider the work of Wei and Jiang as a special case of equation (\ref{VanDerPol}). We will concentrate on the Zero-Hopf bifurcation, and we will show that the added delay in the velocity has now allowed us to fully realize all four generic unfoldings of this singularity (in contrast to only two possible unfoldings in \cite{WuWang}, as mentioned above). We will use methods of normal form reduction on the centre manifold to effectively study the dynamics of the system near the bifurcation. This will therefore allow us to observe the dynamics of the Zero-Hopf bifurcation as it occurs on the centre manifold. 

The paper is organized as follows.  We begin with an analysis of the characteristic equation of equation (\ref{VanDerPol}). 
We then use a standard functional analytic framework \cite{FM1,FM2} to reduce (\ref{VanDerPol}) to a normal form on the three-dimensional center manifold for the Zero-Hopf bifurcation.
We will then use the well-known analysis of the Zero-Hopf bifurcation \cite{GuckHolmes,Kuznetsov} to completely classify the possible local phase diagrams near the singularity in terms of the model parameters in (\ref{VanDerPol}).  Finally, we will present numerical results illustrating the theoretical results of the paper.  Some technical aspects of our work are detailed in the Appendix.

\Section{Linear Analysis}

Linearization of equation (\ref{VanDerPol}) at the trivial solution $x=0$ gives
\begin{equation}
\ddot{x}(t)-\varepsilon \dot{x}(t)+x(t)=a\dot{x}(t-\tau)+b x(t-\tau).
\label{LinVanDerPol}
\end{equation}
Substituting the ansatz $x(t)=e^{\lambda t}$ into (\ref{LinVanDerPol}) gives us the characteristic equation
\begin{equation}
\Delta(\lambda,\tau):=\lambda^2-\varepsilon\lambda+1-(a\lambda+b)e^{-\lambda\tau}=0.
\label{CharacteristicEquation}
\end{equation}
We then get the following

\begin{prop} \label{EigenvalueTheorem}
Suppose that $b = 1$ is satisfied. Then
\begin{itemize}
	\item[(i)] $\lambda = 0$ is a single root to (\ref{CharacteristicEquation}) when $\tau \neq \varepsilon + a$;
	\item[(ii)] $\lambda = 0$ is a double root to (\ref{CharacteristicEquation}) when $\tau = \varepsilon + a$ and $\varepsilon^2 - a^2 \neq 2$;
	\item[(iii)] $\lambda = 0$ is a triple root to (\ref{CharacteristicEquation}) when $\tau = \varepsilon + a$, $\varepsilon^2 - a^2 = 2$ and $\varepsilon \neq \varepsilon_0 \approx 1.632993162$;
	\item[(iv)] $\lambda = 0$ is a quadruple root to (\ref{CharacteristicEquation}) when $\tau = \varepsilon+ a$, $\varepsilon^2 - a^2 = 2$ and $\varepsilon = \varepsilon_0$;
	\item[(v)] If $\tau = \tau_0 \neq \varepsilon + a$ and $\varepsilon^2 - a^2 < 2$, (\ref{CharacteristicEquation}) has a simple zero root and a pair of purely imaginary roots $\lambda = \pm i\omega_0$. Here $\omega_0$ and $\tau_0$ are defined by
	\[
		\omega_0 = \sqrt{2 - \varepsilon^2 + a^2},\ \ \ \ \ \tau_0 = \frac{1}{\omega_0}\arccos\bigg(\frac{1-(1+\varepsilon a)\omega_0^2}{a^2\omega_0^2 + 1}\bigg).
	\] 
\end{itemize}
\end{prop}
\proof
	Immediately we see that $\lambda = 0$ is a root of (\ref{CharacteristicEquation}) if and only if $b = 1$. Substituting $b=1$ into $\Delta(\lambda,\tau)$ and differentiating with respect to $\lambda$, we get
	\begin{equation} \label{FirstDiff} 
		\frac{d\Delta(\lambda,\tau)}{d\lambda} = 2\lambda - \varepsilon - ae^{-\lambda\tau} + a\lambda\tau e^{-\lambda\tau} + \tau e^{-\lambda\tau}.
	\end{equation}
	Hence when $\lambda = 0$, (\ref{FirstDiff}) is 0 if and only if $\tau = \varepsilon + a$, and the conclusion of (i) follows.
	
	By differentiating (\ref{FirstDiff}), we obtain
	\begin{equation} \label{SecondDiff} 
		\frac{d^2\Delta(\lambda,\tau)}{d\lambda^2}\bigg|_{\tau=\varepsilon+a} = 2 - (\varepsilon^2 - a^2)e^{-\lambda(\varepsilon + a)} - a\lambda(\varepsilon + a)^2e^{-\lambda(\varepsilon + a)}.		
	\end{equation}   
	Again setting $\lambda = 0$, (\ref{SecondDiff}) is 0 if and only if $\varepsilon^2 - a^2 = 2$, and hence the conclusion of (ii) follows. 
	
	From (\ref{SecondDiff}), we find that 
	\begin{equation} \label{ThirdDiff} 
		\frac{d^3\Delta(0,\tau)}{d\lambda^3}\bigg|_{\tau=\varepsilon+a, \varepsilon^2 - a^2 = 2, \lambda = 0} = 6\varepsilon - 2\varepsilon^3 \pm 2(\varepsilon^2 - 2)\sqrt{\varepsilon^2-2}.	
	\end{equation}
	Equation (\ref{ThirdDiff}) has only one positive root, denoted $\varepsilon_0$, and the conclusion of (iii) follows.
	
	From (\ref{ThirdDiff}), we get that $d^4\Delta(0,\tau)/d\lambda^4|_{\tau=\varepsilon_0+a, \varepsilon_0^2 + a^2 = 2} \neq 0$ for either choice of $a$, which gives the conclusion of (iv). 
	
	By taking $\lambda = i\omega$ ($\omega > 0$) in (\ref{CharacteristicEquation}), when $b=1$ we get
	\begin{equation} \label{RealImaginary}
		1 - \omega^2 = \cos(\omega\tau) + a\omega \sin(\omega\tau) \ \ \textrm{and} \ \ -\varepsilon\omega = a\omega \cos(\omega\tau) - \sin(\omega\tau).
	\end{equation}
	Squaring and adding the above equations gives $\omega^2[\omega^2 + (\varepsilon^2 - 2 - a^2)] = 0$. Denote the positive, nonzero root of this equation $\omega_0$. Similarly, the expression for $\tau_0$ can be obtained by a simple manipulation on the equations in (\ref{RealImaginary}) to isolate $\cos(\omega\tau)$, and hence the conclusion of (v) follows.   
\hfill\qed

Although it will not be the focus of this paper, we remark that by allowing the forcing function $g$ in (\ref{VanDerPol}) to depend also on delayed velocity leads to the possibility of a quadruple zero eigenvalue.  For our purposes, we will be concentrating on case (v) of the previous Proposition, i.e. the Zero-Hopf singularity.

\Section{Functional Analytic Framework}

In this section, we will adapt the technique and notations of \cite{FM1,FM2} to our particular problem.

We rewrite (\ref{VanDerPol}) as the equivalent first-order system
\begin{equation}
\begin{array}{rcl}
\dot{u}_1(t)&=&u_2(t)\\[0.2in]
\dot{u}_2(t)&=&-u_1(t)-\varepsilon (u_1(t)^2-1) u_2(t)+g(u_2(t-\tau),u_1(t-\tau))
\end{array}
\label{VDPsys}
\end{equation}
Next, we consider for $\tau \geq 0$, the phase space $C = C([-\tau,0];\mathbb{C}^2)$, the space of continuous functions from the interval $[-\tau,0]$ into $\mathbb{C}^2$. 
The system (\ref{VDPsys}) is then viewed as a retarded functional differential equation with parameters
 \begin{equation} \label{FunctionalEquation}
 	\dot{z}(t) = L(\mu)z_t + G(z_t,\mu)
 \end{equation}  
where $\mu \in V$, a neighbourhood of zero in $\mathbb{R}^p$, $z_t(\theta) = z(t+\theta)$, $-\tau \leq \theta \leq 0$, $L: C \to \mathbb{C}^2$ is a bounded linear operator and $G: C \times \mathbb{R}^p \to \mathbb{C}^2$ is a sufficiently smooth function such that $G(0,0) = 0$ and $D_{u}G(0,0) = 0$. 

By denoting $L_0 = L(0)$ the linear homogeneous retarded functional differential equation at $\mu = 0$ can be written
\begin{equation} \label{OperatorL}
	\dot{z}(t) = L_0z_t.
\end{equation}   
From Reisz's Representation Theorem, $L_0$ has the form
\begin{equation} \label{LinearOperator}
	L_0\varphi = \int_{-\tau}^0 d\eta(\theta)\varphi(\theta), \ \ \varphi \in C
\end{equation}
where $\eta(\theta)$, $-\tau \leq \theta \leq 0$, is a $2 \times 2$ matrix whose elements are of bounded variation. Therefore, equation (\ref{OperatorL}) becomes
\begin{equation}
	\dot{z}(t) = \int_{-\tau}^0 d\eta(\theta)z(t+\theta).
\end{equation}  
Let $\mathcal{A}_0$ be the infinitesimal generator of the semigroup of solutions (see \cite{HVL} for details) to equation (\ref{OperatorL}), we have $\mathcal{A}_0\varphi = \dot{\varphi}$ with domain
\begin{equation} \label{InfDomain}
	D(\mathcal{A}_0) = \bigg\{\varphi \in C^1([-\tau,0],\mathbb{C}^2) : \dot{\varphi}(0) = \int_{-\tau}^{0} d\eta(\theta)\varphi(\theta)\bigg\}.
\end{equation}

Define $C^* = C([-\tau,0],\mathbb{C}^{2*})$, where $\mathbb{C}^{2*}$ is the space of row vectors. Then for $\varphi \in C$ and $\psi \in C^{*}$ we define the bilinear form
\begin{equation} \label{BilinearForm}
	\langle \psi,\varphi\rangle = \psi(0)\varphi(0) - \int_{-\tau}^0\int_0^{\theta} \psi(\xi - \theta)d\eta(\theta)\varphi(\xi)d\xi.
\end{equation}
Call $\varLambda$ the set of eigenvalues of the operator $\mathcal{A}_0$ with zero real part, and let $P$ be the generalized eigenspace associated with the eigenvalues in $\varLambda$. 
Then, 
for the purpose of this paper, we will assume that dim$P = m$ is finite.  We denote $P^{*}$ as the generalized eigenspace associated to the eigenvalues in $\varLambda$ of the adjoint of $\mathcal{A}_0$ in $C^{*}$. Now $C$ can be decomposed by $\varLambda$ such that $C = P \oplus Q$, where
\begin{equation} \label{QSpace}
	Q = \{\varphi \in C : \langle \psi, \varphi\rangle = 0 \ {\rm for} \ {\rm all} \ \psi \in P^{*}\}.
\end{equation}
Then there exists bases $\Phi$, $\Psi$ for $P$, $P^{*}$ respectively satisfying $\langle \Psi,\Phi\rangle = I$. Using these bases, we will project solutions from the infinite dimensional space $C$ to the finite dimensional centre manifold $P$.  

As described by Faria and Magalh$\tilde{\rm{a}}$es in $\cite{FM1}$ and $\cite{FM2}$, we consider the phase space $BC$ of functions from $[-\tau,0]$ to $\mathbb{C}^2$ which are uniformly continuous on $[-\tau,0)$ and with a jump discontinuity at 0. Therefore any element in $BC$ can be written in the form $\phi = \varphi + X_0c$, where $X_0$ is the matrix valued function given by  
\begin{equation} \label{X0}
	X_0(\theta) = \left\{
     \begin{array}{lr}
       I, & \theta = 0 \\
       0, & \theta \in [-\tau,0)
     \end{array},
   \right.
\end{equation}
$\varphi \in C$, and $c$ is a vector in $\mathbb{C}^2$. This allows us to consider equation (\ref{FunctionalEquation}) in $BC$ as the abstract ordinary differential equation  
\begin{equation} \label{AbstractODE}
	\frac{d}{dt}u = \mathcal{A}u + X_0F(u,\mu)
\end{equation}
where 
\begin{equation} \label{NewNonlinear}
	F(\varphi,\mu) = (L(\mu) - L_0)\varphi + G(\varphi,\mu) 
\end{equation}
for $\varphi \in C$, $\mu \in V$ and the operator $\mathcal{A}:C^1 \to BC$ is an extension of the infinitesimal generator $\mathcal{A}_0$ and is given by
\begin{equation} \label{OperatorA}
	\mathcal{A}\varphi = \dot{\varphi} + X_0[L_0\varphi - \dot{\varphi}(0)].
\end{equation}
Then $\pi:BC \to P$ is the continuous projection 
\begin{equation} \label{Projection}
	\pi(\varphi + X_0c) = \Phi[\langle\Psi,\varphi\rangle + \Psi(0)c]
\end{equation} 
which allows us to decompose the space $BC$ by $\varLambda$ such that $BC = P\ \oplus$ Ker $\pi$. Since we have that $Q \subsetneq$ Ker $\pi$ we may rewrite equation ($\ref{AbstractODE}$) as
\begin{equation} \label{RewrittenAbstract}
	\begin{split}
	\begin{aligned}	
		&\dot{x} = Jx + \Psi(0)F(\Phi x+y,\mu) \\
		& \frac{d}{dt}y = \mathcal{A}_{Q^1}y + (I - \pi)X_0F(\Phi x + y,\mu)
	\end{aligned}
	\end{split}
\end{equation}  
where $x \in \mathbb{C}^m$, $y \in Q^1 = Q \cap C^1$, $\mathcal{A}_{Q^1}$ is the operator $\mathcal{A}$ restricted to functions in $Q^1$, and the $m \times m$ matrix $J$ satisfies the ordinary differential equations $d\Phi/d\theta = \Phi J$ and $d\Psi/d\theta = -J\Psi$. 

\subsection{Application to the Zero-Hopf Bifurcation in (\ref{VDPsys})}

From Proposition \ref{EigenvalueTheorem}, we have that if $\varepsilon>0$, $\varepsilon^2-a^2<2$, $b=1$ and $\tau=\tau_0\neq\varepsilon+a$, then system (\ref{VDPsys}) has a simple zero eigenvalue and a simple pair of purely imaginary eigenvalues $\pm i\omega_0$.

We write the Taylor expansion of $g$ as 
\begin{equation} \label{Taylor_g}
	\begin{split}
		\begin{aligned}
			g(x,y) = &ay + bx + \frac{1}{2}g_{xx}(0,0)x^2 + g_{xy}(0,0)xy + \frac{1}{2}g_{yy}(0,0)y^2 + \frac{1}{3!}g_{xxx}(0,0)x^3 \\ & + \frac{1}{2}g_{xxy}(0,0)x^2y + \frac{1}{2}g_{yyx}(0,0)xy^2 + 			\frac{1}{3!}g_{yyy}(0,0)y^3 + \mathcal{O} (|x,y|^4).
		\end{aligned}
	\end{split}
\end{equation}
Then, after the scaling $t \rightarrow t/\tau$, system (\ref{VDPsys}) becomes
\begin{equation} \label{ScaledSystem}
    		\begin{array}{lr}
       			\ \dot{u_1} = \tau u_2,\\
			\begin{split}
			\begin{aligned}
       				\dot{u_2} = &-\tau u_1 - \varepsilon \tau (u_1^2 - 1)u_2 + \tau \bigg[au_2(t - 1) + bu_1(t - 1) + \frac{1}{2}g_{u_1u_1}(0,0)u_1^2(t - 1) \\  &+ g_{u_1u_2}(0,0)u_1(t - 1)u_2(t - 1) 					+ \frac{1}{2}g_{u_2u_2}(0,0)u_2^2(t - 1) + \frac{1}{6}g_{u_1u_1u_1}(0,0)u_1^3(t - 1) \\ & + \frac{1}{2}g_{u_1u_1u_2}(0,0)u_1^2(t - 1)u_2(t - 1) + \frac{1}{2}g_{u_2u_2u_1}					(0,0)u_1(t - 1)u_2^2(t - 1) \\ & + \frac{1}{6}g_{u_2u_2u_2}(0,0)u_2^3(t - 1)\bigg] + \mathcal{O} (|(u_1,u_2)|^4).
     			\end{aligned}
			\end{split}
		\end{array}
\end{equation}
If we let $\mu_1 = b -1$ and $\mu_2 = \tau - \tau_0$, then $\mu_1$ and $\mu_2$ become bifurcation parameters and system ($\ref{ScaledSystem}$) now becomes
\begin{equation} \label{BifurcationSystem}
    		\begin{array}{lr}
       			\ \dot{u_1} = \tau_0 u_2 + \mu_2 u_2,\\
			\begin{split}
			\begin{aligned}
       				\dot{u_2} = &-\tau_0 u_1 + \tau_0 u_1(t - 1) + \tau_0 a u_2(t - 1) + \varepsilon \tau_0 u_2 - \mu_2 u_1\\ &+ \varepsilon \mu_2 u_2 + a\mu_2 u_2(t - 1) + \tau_0 \mu_1 u_1(t - 				1) + \mu_2 u_1(t - 1) - \varepsilon \tau_0 u_1^2 u_2- \varepsilon \mu_2 u_1^2 u_2\\ &+ \mu_1 \mu_2 u_1(t - 1) + \frac{1}{2}\tau_0g_{u_1u_1}(0,0)u_1^2(t - 1) + 							\tau_0g_{u_1u_2}(0,0)u_1(t - 1)u_2(t - \tau)\\ &+ \frac{1}{2}\tau_0g_{u_2u_2}(0,0)u_2^2(t - 1) + \frac{1}{2}g_{u_1u_1}(0,0)\mu_2u_1^2(t - 1)\\ &+ g_{u_1u_2}(0,0)\mu_2u_1(t - 1)u_2(t - 1) + \frac{1}{2}g_{u_2u_2}(0,0)\mu_2u_2^2(t - 1)\\ &+ \frac{1}{6}\tau_0g_{u_1u_1u_1}(0,0)u_1^3(t - 1) + \frac{1}{2}\tau_0g_{u_1u_1u_2}(0,0)u_1^2(t - 1)u_2(t - 1)\\ &+ \frac{1}{2}\tau_0g_{u_2u_2u_1}(0,0)u_1(t - 1)u_2^2(t - 1) + \frac{1}{6}\tau_0g_{u_2u_2u_2}(0,0)u_2^3(t - 1) +h.o.t., 
	     		\end{aligned}
			\end{split}
		\end{array}
\end{equation}
where $h.o.t.$ stands for ``higher order terms'' and in this case represents the terms which are of order $\mathcal{O}(|(u_1,u_2)|^4)$ uniformly in $\mu_1$ and $\mu_2$.

Linearizing (\ref{BifurcationSystem}) around $(u_1,u_2,\mu_1,\mu_2) = (0,0,0,0)$ we get
\begin{equation} \label{LinearBifurcationSystem}
    		\begin{array}{lr}
       			\ \dot{u_1} = \tau_0 u_2, \\
			\begin{split}
			\begin{aligned}
       				\dot{u_2} = &-\tau_0 u_1 + \varepsilon \tau_0 u_2  + \tau_0 u_1(t - 1) + \tau_0 a u_2(t - 1). 
	     		\end{aligned}
			\end{split}
		\end{array}
\end{equation}
By letting 
\begin{equation} \label{EtaEq}
	d\eta(\theta) = \mathbb{A}\delta(\theta) + \mathbb{B}\delta(\theta + 1)	
\end{equation}
such that 
\begin{equation} \label{AandB}
	\mathbb{A} = \begin{pmatrix} 0 & \tau_0 \\ -\tau_0 & \varepsilon\tau_0 \end{pmatrix}, \ \ \ \ \ \mathbb{B} = \begin{pmatrix} 0 & 0 \\ \tau_0 & a\tau_0 \end{pmatrix},
\end{equation}
we can define the operator 
\begin{equation}
	L\varphi = \int_{-1}^0 d\eta(\theta)\varphi(\theta), \ \ \ \forall \varphi \in C = C([-1,0],\mathbb{C}^2). 
\end{equation} 
Define $X = (u_1,u_2)^T$. Then (\ref{BifurcationSystem}) can be rewritten in the following form,
\begin{equation} \label{AbstractBifurcationSystem}
		\dot{X}(t) = LX_t + F(X_t,\mu),	
\end{equation}  
where $F(X_t,\mu) = (F^1,F^2)^T$ is given by
\begin{equation} 
	\begin{array}{lr}
       		\ F^1 =  \mu_2 u_2,\\
			\begin{split}
			\begin{aligned}
       				F^2 = &\mu_2 u_1 + \varepsilon \mu_2 u_2 + a\mu_2 u_2(t - 1) + \tau_0 \mu_1 u_1(t - 1) + \mu_2 u_1(t - 1)\\ &- \varepsilon \tau_0 u_1^2 u_2 - \varepsilon \mu_2 u_1^2 				u_2 + \mu_1 \mu_2 u_1(t - 1) + \frac{1}{2}\tau_0g_{u_1u_1}(0,0)u_1^2(t - 1)\\ &+ \tau_0g_{u_1u_2}(0,0)u_1(t - 1)u_2(t - 1) + \frac{1}{2}\tau_0g_{u_2u_2}(0,0)u_2^2(t - 1) 				+\\ &\frac{1}{2}g_{u_1u_1}(0,0)\mu_2u_1^2(t - 1) + g_{u_1u_2}(0,0)\mu_2u_1(t - 1)u_2(t - 1) + \frac{1}{2}g_{u_2u_2}(0,0)\mu_2u_2^2(t - 1)\\ &+ \frac{1}								{6}\tau_0g_{u_1u_1u_1}(0,0)u_1^3(t - 1) + \frac{1}{2}\tau_0g_{u_1u_1u_2}(0,0)u_1^2(t - 1)u_2(t - 1)\\ &+ \frac{1}{2}\tau_0g_{u_2u_2u_1}(0,0)u_1(t - 1)u_2^2(t - 1) + 					\frac{1}{6}\tau_0g_{u_2u_2u_2}(0,0)u_2^3(t - 1) +h.o.t.
	     		\end{aligned}
			\end{split}
		\end{array}	
\end{equation}
Using the definition of $d\eta$ above we can simplify the bilinear form in (\ref{BilinearForm}) to get 
\begin{equation} \label{MyBilinearForm}
	\langle \psi, \varphi \rangle = \psi(0)\varphi(0) + \int_{-1}^{0} \psi(\xi + 1)\mathbb{B}\varphi(\xi)d\xi.
\end{equation}
The infinitesimal generator $\mathcal{A}: C^1([-1,0],\mathbb{C}^2) \to BC$ is defined as
\begin{equation} \label{InfGenerator}
	\mathcal{A}\varphi(\theta) = \dot{\varphi} + X_0[L\varphi - \dot{\varphi}(0)] = \left\{
    		\begin{array}{ll}
       			\dot{\varphi} & \mathrm{if}\ -1 \leq \theta < 0, \\ 
			{\int_{-1}^0} d\eta(t)\varphi(t), & \mathrm{if}\  \theta = 0,
		\end{array}
   	\right.	
\end{equation}
with adjoint $\mathcal{A}^*: C^{1}([0,1],\mathbb{C}^{2*}) \to BC$ given by
\begin{equation} \label{AdjointInfGenerator}
	\mathcal{A}^*\psi(s) = -\dot{\psi} + X_0[L^*\psi + \dot{\psi}(0)] = \left\{
    		\begin{array}{ll}
       			\dot{\psi} & \mathrm{if} \ 0 \leq s < 1, \\ 
			{\int_{0}^1} \psi(-t)d\eta(t), & {\rm if} \ s = 0.
		\end{array}
   	\right.
\end{equation}
Then $\mathcal{A}$ and $\mathcal{A}^*$ have eigenvalues $0, i\omega_0\tau_0$ and $-i\omega_0\tau_0$ for which we must now compute their respective eigenvectors.

\begin{prop} \label{Eigenvectors} 
	There exists bases $\Phi = (\varphi_1,\bar{\varphi}_1,\varphi_2)$ and $\Psi = (\bar{\psi}_1,\psi_1,\psi_2)^T$ of the centre space $P$ and its adjoint $P^*$ such that $\dot{\Phi} = \Phi J$, $\dot{\Psi} = J\Psi$ and $\langle \Psi,\Phi \rangle = I$, where
\[
\begin{split}
\begin{aligned}
	&\varphi_1(\theta) = (1, i\omega_0)^Te^{i\omega_0\tau_0\theta}, \ \ \ \ \ \varphi_2(\theta) = (1,0)^T \\
	&\psi_1(s) = D(1, \sigma)e^{i\omega_0\tau_0s}, \ \ \ \ \ \psi_2(s) = D_1(\varepsilon + a, -1) \\
	&J = diag(i\omega_0\tau_0, -i\omega_0\tau_0, 0), \ \ \ \ \ \sigma = \frac{i\omega_0}{1 - e^{i\omega_0\tau_0}},  \\
	&D = \frac{1}{1 - i\sigma\omega_0 + \tau_0\sigma e^{i\omega_0\tau_0}(1 - ia\omega_0)}, \ \ \ \ \ D_1 = \frac{1}{\varepsilon + a - \tau_0}.
\end{aligned}
\end{split}
\]    
\end{prop}
\proof
	Suppose $\varphi_1(\theta) = (1, \rho)^Te^{i\omega_0\tau_0\theta}$ is an eigenvector of $\mathcal{A}$ corresponding to $i\omega_0\tau_0$. Then $\mathcal{A}\varphi_1(\theta) = i\omega_0\tau_0\varphi_1(\theta)$. From the definition of $\mathcal{A}$ we get \\
\[
	\mathbb{A}\varphi_1(0) + \mathbb{B}\varphi_1(-1) = i\omega_0\tau_0\varphi_1(0)
\] \\	
which gives us that \\
\[
	\tau_0
	\begin{pmatrix}	
		-i\omega_0 & 1 \\
		-1 + e^{-i\omega_0\tau_0} & \varepsilon - ae^{-i\omega_0\tau_0} - i\omega_0 	
	\end{pmatrix}
	\begin{pmatrix}
		1 \\
		\rho
	\end{pmatrix}
	=
	\begin{pmatrix}
		0 \\
		0
	\end{pmatrix},
\] \\
which gives $\rho = i\omega_0$. Similarly suppose $\varphi_2(\theta)$ is the eigenvector of $\mathcal{A}$ corresponding to 0. From the definition of $\mathcal{A}$, we have that $\varphi_2(\theta)$ is a constant vector $(\alpha,\beta)^T$ such that \\
\[
	(\mathbb{A} + \mathbb{B})
	\begin{pmatrix}
		\alpha \\
		\beta
	\end{pmatrix}
	= 0
\] \\
from which we obtain $\varphi_2(\theta) = (1,0)^T$. Let $\Phi = (\varphi_1(\theta), \bar{\varphi_1}(\theta),\varphi_2(\theta))$. Then it is easy to check that $\dot{\Phi} = \Phi J$.

Now suppose $\psi_1(s) = D(1, \sigma)e^{i\omega_0\tau_0s}$ is an eigenvector of $\mathcal{A}^*$ corresponding to $-i\omega_0\tau_0$. Then $\mathcal{A}^*\psi_1(s) = -i\omega_0\tau_0\psi_1(s)$. Then from the definition of $\mathcal{A}^*$ we have \\
 \[
 	\psi_1(0)\mathbb{A} + \psi_1(1)\mathbb{B} = -i\omega_0\tau_0\psi_1(0)
 \] \\
 or equivalently \\
 \[
 	\tau_0
	\begin{pmatrix}
		1, & \sigma
	\end{pmatrix}
	\begin{pmatrix}
		i\omega_0 & 1 \\
		-1 + e^{i\omega_0\tau_0} & \varepsilon + ae^{i\omega_0\tau_0} + i\omega_0
	\end{pmatrix} = (0,0)
 \] \\
 which gives us \\
 \[
 	\sigma = \frac{i\omega_0}{1 - e^{i\omega_0\tau_0}}.
 \] \\
 Suppose that $\psi_2(s)$ is the eigenvector of $\mathcal{A}^*$ corresponding to 0. From the definition of $\mathcal{A}^*$, we have that $\psi_2(s)$ is a constant vector $(\alpha, \beta)$ such that \\
 \[
 	(\alpha, \beta)(\mathbb{A} + \mathbb{B}) = 0
 \] \\
 which gives $\psi_2(s) = D_1(\varepsilon + a, -1)$. Let $\Psi = (\bar{\psi_1}(s), \psi_1(s),\psi_2(s))$. Then we have $\dot{\Psi} = -J\Psi$.
 
 It can be checked that $\langle \psi_1,\varphi_1\rangle = 0$, $\langle \psi_2, \varphi_1\rangle = 0$ and $\langle \psi_1, \varphi_2\rangle = 0$. In order to assure $\langle \bar{\psi_1},\varphi_1 \rangle = 1$, $\langle \psi_2,\varphi_2 \rangle = 1$, we must determine the factors $D, D_1$. Since \\
 \[
 \begin{split}
 \begin{aligned}
  	\langle \bar{\psi_1},\varphi_1 \rangle &= \bar{D} \bigg[(1,\bar{\sigma})
	\begin{pmatrix}
		1 \\
		i\omega_0
	\end{pmatrix} 
	+ \int_{-1}^0 (1, \bar{\sigma})e^{-i\omega_0\tau_0(\xi + 1)}\mathbb{B}
	\begin{pmatrix}
		1 \\
		i\omega_0
	\end{pmatrix}
	e^{i\omega_0\tau_0\xi}d\xi\bigg], \\
	&= \bar{D}\big[1 + i\bar{\sigma}\omega_0 + \tau_0\bar{\sigma}e^{-i\omega_0\tau_0}(1 + ia\omega_0)\big]  \\
	\langle \psi_2, \varphi_2 \rangle &= D_1 \bigg[(\varepsilon + a, -1)
	\begin{pmatrix}
		1 \\
		0
	\end{pmatrix}
	+ \int_{-1}^0 (\varepsilon + a, -1)\mathbb{B}
	\begin{pmatrix}
		1 \\
		0
	\end{pmatrix}
	d\xi \bigg] \\
	&= D_1\big[\varepsilon + a - \tau_0\big],
\end{aligned}
\end{split}
 \]  \\
 we get the desired results, finishing the proof of the lemma.
\hfill\qed	

Now let $u = \Phi x + y$ where $x \in \mathbb{C}^3$ and $y \in Q^1 = \{\varphi \in Q : \dot{\varphi} \in C\}$, then we write $u$ as
\begin{equation}
	\begin{split}
	\begin{aligned}
		u_1(\theta) &= e^{i\omega_0\tau_0\theta}x_1 + e^{-i\omega_0\tau_0\theta}x_2 + x_3 + y_1(\theta) \\
		u_2(\theta) &= i\omega_0e^{i\omega_0\tau_0\theta}x_1 - i\omega_0e^{-i\omega_0\tau_0\theta}x_2 + y_2(\theta),
	\end{aligned}
	\end{split}
\end{equation} 
with $x_2=\overline{x_1}$.
Let
\begin{equation}
	\Psi (0) =
	\begin{pmatrix}
		\psi_{11} & \psi_{12} \\
		\psi_{21} & \psi_{22} \\
		\psi_{31} & \psi_{32} 
	\end{pmatrix}
	=
	\begin{pmatrix}
		\bar{D} & \bar{D}\bar{\sigma} \\
		D & D\sigma \\
		D_1(\varepsilon + a) & -D_1
	\end{pmatrix}.
\end{equation}
Then System ($\ref{BifurcationSystem}$) is decomposed as
\begin{equation} \label{DecomposedSystem}  
\left\{
	\begin{split}
	\begin{aligned}
		&\dot{x} = Jx + \Psi(0)F(\Phi x + y, \mu) \\
		&\dot{y} = A_{Q^1}y + (I - \pi)X_0F(\Phi x + y, \mu)
	\end{aligned}
	\end{split}
\right.
\end{equation}
and upon using Taylor's Theorem we obtain 
\begin{equation} \label{ExpandedDecomposedSystem}
	 \left\{
	\begin{split}
	\begin{aligned}
		&\dot{x} = Jx + f^1_2(x,y,\mu) + f_3^1(x,y,\mu) + h.o.t., \\
		&\dot{y} = A_{Q^1}y + f^2_2(x,y,\mu) + f_3^2(x,y,\mu) + h.o.t.,
	\end{aligned}
	\end{split}
	\right.	
\end{equation}
where for $j=2,3$
\begin{equation} \label{TaylorTerms1}
\begin{split}
 \begin{aligned}
 	&f_j^1(x,y,\mu) = 
	\begin{pmatrix}
		\psi_{11}F_j^1(\Phi x + y, \mu) + \psi_{12}F^2_j(\Phi x + y, \mu) \\
		\psi_{21}F_j^1(\Phi x + y, \mu) + \psi_{22}F^2_j(\Phi x + y, \mu) \\
		\psi_{31}F_j^1(\Phi x + y, \mu) + \psi_{32}F^2_j(\Phi x + y, \mu) 
	\end{pmatrix}, \\
	&f_2^2(x,y,\mu) = (I - \pi)X_0
	\begin{pmatrix}
		F^1_j(\Phi x + y, \mu) \\
		F^2_j(\Phi x + y, \mu)
	\end{pmatrix}.
 \end{aligned}
 \end{split}
\end{equation}

\subsection{Faria and Magalh$\tilde{\mbox{\rm a}}$es Normal Form}

We are now ready to state our main result of this paper.

\begin{thm}
Near the triple-zero nilpotent bifurcation point in (\ref{VDPsys}), there exists a 3-dimensional invariant center manifold on which the local dynamics of (\ref{VDPsys}) reduce to the 3-dimensional ordinary differential equation
\begin{equation} \label{CentreODE}
\begin{split}
\begin{aligned}
       		\dot{x_1} = &i\omega_0x_1 + (a_{11}\mu_1 + a_{12}\mu_2)x_1 + a_{13}x_1x_3 + (b_{11} + c_{11} + d_{11} - e_{11} - m_{11} + n_{11})x_1^2x_2\\ &+ (b_{12} + c_{12} + d_{12} - e_{12} - m_{12} + n_{12})x_1x_3^2 + h.o.t., \\
		\dot{x_2} = &-i\omega_0x_2 + (\bar{a}_{11}\mu_1 + \bar{a}_{12}\mu_2)x_2 + \bar{a}_{13}x_2x_3 + (\bar{b}_{11} + \bar{c}_{11} + \bar{d}_{11} - \bar{e}_{11} -\bar{m}_{11} + \bar{n}_{11} )x_1x_2^2\\ &+ (\bar{b}_{12} + \bar{c}_{12} + \bar{d}_{12} - \bar{e}_{12} - \bar{m}_{12} + \bar{n}_{12})x_2x_3^2 + h.o.t.,  \\ 
		\dot{x_3} = &(a_{21}\mu_1 + a_{22}\mu_2)x_3 + a_{23}x_3^2 + (b_{21} + c_{21} + d_{21} - e_{21} - m_{21} + n_{21})x_1x_2x_3\\ &+ (b_{22} + c_{22} + d_{22} - e_{22} - m_{22} + n_{22})x_3^3 + h.o.t
\end{aligned}
\end{split}
\end{equation}
where
$x_1\in\mathbb{C}$, $x_2=\overline{x}_1$, and $x_3\in\mathbb{R}$, 
where h.o.t. refers to 
terms which are $\mathcal{O}(|(x_1,x_2,x_3)|^4)$ uniformly in $\mu$,
and where the coefficients $a_{ij}$, $b_{ij}, \ldots , m_{ij}, n_{ij}$ are given in terms of the various Maclaurin coefficients of the function $g$ in (\ref{VDPsys}) in the Appendix.\label{MainThm}
\end{thm}
\proof
The proof hinges on performing, at successively higher orders, near identity changes of variables of the form
\[
(x,y)=(\hat{x},\hat{y})+(U_1(\hat{x}),U_2(\hat{x}))
\]
on (\ref{ExpandedDecomposedSystem})
as per \cite{FM1,FM2}, and setting $y=0$.  Although the computational details are lengthy and tedious, they are standard and straightforward.  The details are given in the Appendix.
\hfill\qed

\Section{Reduction to Polar Coordinates}

We can convert the ODE (\ref{CentreODE}) from complex variables to real variables by introducing the change of variables $x_1 = w_1 - iw_2$, $x_2 = w_1 + iw_2$, $x_3 = z$. Upon completing this change of variables we may now convert to cylindrical coordinates by letting $w_1 = r\cos\xi$, $w_2 = r\sin\xi$. System ($\ref{CentreODE}$) now becomes 
\begin{equation} \label{CylindricalODE}
	\left\{
    	\begin{array}{ll}
       		\dot{r} = \alpha_1(\mu)r + \beta_{11}rz + \beta_{30}r^3 + \beta_{12}rz^2 + h.o.t., \\
		\dot{z} = \alpha_2(\mu)z + \gamma_{20}r^2 + \gamma_{02}z^2 + \gamma_{03}z^3 + h.o.t., \\
		\dot{\xi} = -\omega_0 + (\mathrm{Im}[a_{11}]\mu_1 + \mathrm{Im}[a_{12}]\mu_2)r + h.o.t.,
	\end{array}
\right.	
\end{equation}   
where
\begin{equation}  \label{CylindricalCoefficients}
	\begin{split}
\begin{aligned}
	&\alpha_1(\mu) = \mathrm{Re}[a_{11}]\mu_1 + \mathrm{Re}[a_{12}]\mu_2, \ \alpha_2(\mu) = a_{21}\mu_1 + a_{22}\mu_2, \ \beta_{11} = \mathrm{Re}[a_{13}], \\ 
	&\beta_{30} = \mathrm{Re}[b_{11} + c_{11} + d_{11} - e_{11} - m_{11} + n_{11}], \ \beta_{12} = \mathrm{Re}[b_{12} + c_{12} + d_{12} - e_{12} - m_{12} + n_{12}], \\ 
	&\gamma_{20} = a_{23}, \ \ \ \gamma_{02} = a_{24}, \ \gamma_{21} = b_{21} + c_{21} + d_{21} - e_{21} - m_{21} + n_{21},\\ &\gamma_{03} = b_{22} + c_{22} + d_{22} - e_{22} -m_{22} + n_{22}.	
\end{aligned}
	\end{split}
\end{equation}

In the sequel, we will ignore the uncoupled $\dot{\xi}$ equation in (\ref{CylindricalODE}) and concentrate on the amplitude equations
\begin{equation} \label{NormalForm2}
	\begin{split}
	\begin{aligned}
		\dot{r} = &\alpha_1(\mu)r + \beta_{11}rz + \beta_{30}r^3 + \beta_{12}rz^2 + h.o.t. \\
		\dot{z} = &\alpha_2(\mu)z + \gamma_{20}r^2 + \gamma_{02}z^2 + \gamma_{03}z^3 + h.o.t.
	\end{aligned}
	\end{split}
\end{equation}

Equilibrium points with $r=0$ for (\ref{NormalForm2}) correspond to equilibrium solutions for (\ref{CylindricalODE}) and (\ref{VDPsys}), equilibrium points with $r>0$ for (\ref{NormalForm2}) correspond to periodic solutions for (\ref{CylindricalODE}) and (\ref{VDPsys}), and limit cycles for (\ref{NormalForm2}) correspond to invariant tori for (\ref{CylindricalODE}) and (\ref{VDPsys}).

We consider a change of variables $z \rightarrow z + \delta$, where $\delta$ is yet to be defined explicitly. This gives 
\begin{equation} \label{ShiftedCylindrical}
	\begin{split}
	\begin{aligned}
\dot{r}=&\Omega_1(\mu,\delta)r+\Omega_2(\mu,\delta)rz + h.o.t\\
\dot{z}=&\delta\Omega_3(\mu,\delta)+\Omega_4(\mu,\delta)z+\Omega_5(\mu,\delta)r^2+\Omega_6(\mu,\delta)z^2 + h.o.t.,
	\end{aligned}
	\end{split}
\end{equation}
where the $\Omega_j$ are smooth functions, and have leading order expansions given by
\[
\begin{array}{lll}
\Omega_1(\mu,\delta)=\alpha_1(\mu)+\beta_{11}\delta+h.o.t.&\Omega_2(\mu,\delta)=\beta_{11}+h.o.t.&\Omega_3(\mu,\delta)=\alpha_2(\mu)+\gamma_{02}\delta+h.o.t.\\
\Omega_4(\mu,\delta)=\alpha_2(\mu)+2\gamma_{02}\delta+h.o.t.&
\Omega_5(\mu,\delta)=\gamma_{20}+h.o.t.&
\Omega_6(\mu,\delta)=\gamma_{02}+h.o.t.
\end{array}
\]
We will assume the following generic conditions:
\begin{hyp}
The function $g$ in (\ref{VDPsys}) is such that the coefficients $\beta_{11}$, $\gamma_{20}$ and $\gamma_{02}$ (via the formulae given in the appendix and in (\ref{CylindricalCoefficients})) are non-zero.  That is, we have
\[
\mbox{\rm Re}[a_{13}]\neq 0,\,\,\,a_{23}\neq 0,\,\,\,a_{24}\neq 0.
\]
\end{hyp}
We may then use the implicit function theorem to solve uniquely the equation $\Omega_4(\mu,\delta)=0$ for $\delta=\delta(\mu)$ near $\mu=0$, $\delta=0$
\begin{equation}
\label{Delta}
\delta(\mu)=-\frac{1}{2\gamma_{02}}\alpha_2(\mu)+\mathcal{O}(|\mu|^2)=-\frac{a_{21}}{2a_{24}}\mu_1+\mathcal{O}(|\mu|^2)
\end{equation}
Equation (\ref{ShiftedCylindrical}) then becomes
\begin{equation} \label{ShiftedCylindrical2}
	\begin{split}
	\begin{aligned}
\dot{r}=&\eta_1(\mu)r+\mbox{\rm Re}[a_{13}]rz + h.o.t\\
\dot{z}=&\eta_2(\mu)+a_{23}r^2+a_{24}z^2 + h.o.t.,
	\end{aligned}
	\end{split}
\end{equation}
where
\begin{equation}
\label{etadef}
\begin{array}{rcl}
\eta_1(\mu)&=&{\displaystyle\left(\mbox{\rm Re}[a_{11}]-\frac{a_{21}}{2a_{24}}\mbox{\rm Re}[a_{13}]\right)\mu_1+\mbox{\rm Re}[a_{12}]\mu_2}\\
\eta_2(\mu)&=&{\displaystyle -\frac{a_{21}^2}{4a_{24}}\mu_1^2}
\end{array}
\end{equation}
Now taking $r \rightarrow -\sqrt{|a_{23}a_{24}|}r$ and $z \rightarrow -a_{24}z$ and truncating to quadratic order, we get 
\begin{equation} \label{NormalForm3}
	\begin{split}
	\begin{aligned}
		\dot{r} = &\chi_1r + Arz, \\
		\dot{z} = &\chi_2 + Br^2 - z^2,
	\end{aligned}
	\end{split}
\end{equation}
where
\begin{equation}
\begin{array}{cc}
		\chi_1 = \eta_1(\mu), & {\displaystyle \chi_2 = -a_{24}\eta_2(\mu)=\frac{a_{21}^2}{4}\mu_1^2,}\\
		{\displaystyle A = -\frac{\mbox{\rm Re}[a_{13}]}{a_{24}},}  &B = \pm 1=-\mbox{\rm sgn}(a_{23}a_{24}).
	\end{array}
	\label{ABdefs}
\end{equation}
\begin{rmk}
We make the following remarks concerning (\ref{NormalForm3}) and the relation between our paper and the paper of Wu and Wang \cite{WuWang}
\begin{enumerate}
\item[(a)]
System (\ref{NormalForm3}) is the standard quadratic order normal form for the Zero-Hopf bifurcation, as can be found in \cite{GuckHolmes,Kuznetsov}.  Note that
if the function $g$ in (\ref{VDPsys}) depends only on $u_1(t-\tau)$, as was the case in \cite{WuWang},
then from the formulae (\ref{aValues}) in the appendix, we get that
\[
B=-\mbox{\rm sgn}(a_{23}a_{24})=\mbox{\rm sgn}(-\frac{1}{2}a_{21}^2g_{u_1,u_1}(0,0)^2) = -1.
\]
Thus, the cases where $B=+1$ in the standard unfolding of the Zero-Hopf bifurcation are not realizable if the forcing function $g$ in (\ref{VDPsys}) depends only on a delayed position term, as was assumed in \cite{WuWang}.  By allowing $g$ to also depend on delayed velocity as we have done in this paper, we can attain the unfoldings with $B=+1$ as well as the unfoldings with $B=-1$.  
\item[(b)]
We also note that since
\[
\chi_2=\frac{a_{21}^2}{4}\mu_1^2
\]
then the phase diagrams with $\chi_2<0$ in the unfoldings of the Zero-Hopf bifurcation are not attainable for our problem.  This is due to the fact that (\ref{VDPsys}) always possesses a trivial equilibrium solution, and any bifurcation from this trivial solution to another equilibrium solution with $u_2=0$ in (\ref{BifurcationSystem}) is a transcritical bifurcation; whereas (\ref{NormalForm3}) exhibits a saddle-node bifurcation on the axis $r=0$ when $\chi_2$ crosses from positive to negative.  This issue is independent of whether or not the function $g$ in (\ref{VDPsys}) depends on only delayed position, or on delayed position and delayed velocity.  However, there is an error on page 2597 of \cite{WuWang} in the formula for $\chi_2$, where they claim that the leading order dependence of $\chi_2$ on $\mu_1$ is of cubic order.
\end{enumerate}
\label{rmk42}
\end{rmk}

Following the analysis done in \cite{GuckHolmes,Kuznetsov}, we know that system (\ref{NormalForm3}) undergoes a pitchfork bifurcation from an $r=0$ steady-state to an $r>0$ steady-state when the parameters $(\chi_1,\chi_2)$ are such that $\chi_2=\chi_1^2/A^2$.  Translating back to our parameters, we get
\begin{prop}
Near $(\mu_1,\mu_2)=(0,0)$, system (\ref{BifurcationSystem}) undergoes a Hopf bifurcation from a steady-state solution to a periodic solution when $(\mu_1,\mu_2)$ is on one of the following bifurcation curves
\[
\mbox{\rm HB1} : \mbox{\rm Re}[a_{11}]\mu_1+\mbox{\rm Re}[a_{12}]\mu_2+\mathcal{O}(|\mu|^2)=0,
\]
\[
\mbox{\rm HB2} : (a_{21}\mbox{\rm Re}[a_{13}]-a_{24}\mbox{\rm Re}[a_{11}])\mu_1-a_{24}\mbox{\rm Re}[a_{12}]\mu_2+\mathcal{O}(|\mu|^2)=0.
\]
\end{prop}
We also know from \cite{GuckHolmes,Kuznetsov} that in the case where $B=-1$ and $A>0$ in (\ref{NormalForm3}) then on the semi-axis $\chi_1=0$, $\chi_2>0$, an $r>0$ steady-state of (\ref{NormalForm3}) undergoes a Hopf bifurcation to a limit-cycle.  Furthermore, this limit-cycle grows in amplitude (as $\chi_1$ moves away from the $\chi_1=0$ axis), and eventually disappears in a global heteroclinic bifurcation on a curve ${\displaystyle\chi_1=-\frac{A}{3A+2}\chi_2+\mathcal{O}(|\chi_2|^\frac{3}{2})}$.
For our parameters, we get the following correspondence
\begin{prop}
Suppose system (\ref{BifurcationSystem}) is such that the coefficients in (\ref{aValues}) satisfy $a_{23}a_{24}>0$ and $\mbox{\rm Re}[a_{13}]/a_{24}<0$.  Then near $(\mu_1,\mu_2)=(0,0)$, system (\ref{BifurcationSystem}) undergoes a torus bifurcation from a periodic solution to an invariant two-torus when $(\mu_1,\mu_2)$ is on the bifurcation curve
\[
\mbox{\rm TB} : {\displaystyle\left(\mbox{\rm Re}[a_{11}]-\frac{a_{21}}{2a_{24}}\mbox{\rm Re}[a_{13}]\right)\mu_1+\mbox{\rm Re}[a_{12}]\mu_2}+\mathcal{O}(|\mu|^2)=0.
\]
In terms of the original parameters of (\ref{BifurcationSystem}), the heteroclinic bifurcation curve for the normal form (\ref{NormalForm3}) is given by
\begin{equation}
\label{hetbifcurve}
\mbox{\rm HET} : {\displaystyle\left(\mbox{\rm Re}[a_{11}]-\frac{a_{21}}{2a_{24}}\mbox{\rm Re}[a_{13}]\right)\mu_1+\mbox{\rm Re}[a_{12}]\mu_2}=\frac{a_{21}^2\mbox{\rm Re}[a_{13}]}{8a_{24}+12\mbox{\rm Re}[a_{13}]}\mu_1^2+\mathcal{O}(\mu_1^3).
\end{equation}
Because higher order terms break the normal form symmetry, the global heteroclinic bifurcation occuring in the normal form (\ref{NormalForm3}) on the curve (\ref{hetbifcurve}) may lead to chaotic dynamics in the system (\ref{BifurcationSystem}) near the torus for parameter values near the curve
(\ref{hetbifcurve}).
\end{prop}

\Section{Numerical Examples}

In this section, we provide results of numerical simulations on the system
\begin{equation}
\ddot{x}(t) + \varepsilon(x^2(t) - 1)\dot{x}(t) + x(t) = g(\dot{x}(t-\tau),x(t-\tau))
\label{simsys}
\end{equation}
for various values of $\varepsilon$ and various functions $g$ for parameter values near the Zero-Hopf point, with the goal of illustrating some of the cases in the unfolding space of this bifurcation which are achieved in (\ref{simsys}).  In all cases, we have used {\em Matlab}; in particular, the {\em dde23} package to integrate the system (\ref{simsys}), and the {\em supsmu} package to represent the solutions graphically.

\subsection{Case I: $A<0$ and $B=1$}

We consider the function
\begin{equation}
g(\dot{x}(t-\tau),x(t-\tau))=x(t-\tau)+0.1\, \dot{x}(t-\tau)-0.2\, x(t-\tau)^2+0.2\, x(t-\tau)\dot{x}(t-\tau)+0.2\, \dot{x}(t-\tau)^2
\label{gnum1}
\end{equation}
and
$
\varepsilon=0.3
$
in (\ref{simsys}).  In this case, a lengthy but straightforward computation of the coefficients $A$ and $B$ as per formulae (\ref{ABdefs}) and (\ref{aValues}) yields
\[
A\approx -3.358,\,\,\,\,\,\,B=1.
\]
In this case, system (\ref{simsys}) undergoes a Zero-Hopf bifurcation
with critical values
\[
\omega_0\approx 1.386,\,\,\,\,\,\tau_0\approx 2.060
\]
as given in Proposition
\ref{EigenvalueTheorem} (v).

We thus introduce the unfolding parameters $\mu_1$ and $\mu_2$ and consider the following unfolding of (\ref{simsys})
\begin{equation}
\begin{array}{l}
\ddot{x}(t) + \varepsilon(x^2(t) - 1)\dot{x}(t) + x(t) = \\
\\
(1+\mu_1)x(t-\tau)+0.1\, \dot{x}(t-\tau)-0.2\, x(t-\tau)^2+0.2\, x(t-\tau)\dot{x}(t-\tau)+0.2\, \dot{x}(t-\tau)^2
\end{array}
\label{simsys2}
\end{equation}
where
\[
\tau=\tau_0+\mu_2.
\]

\begin{figure}[htbp]
\begin{center}
\includegraphics[width=6in]{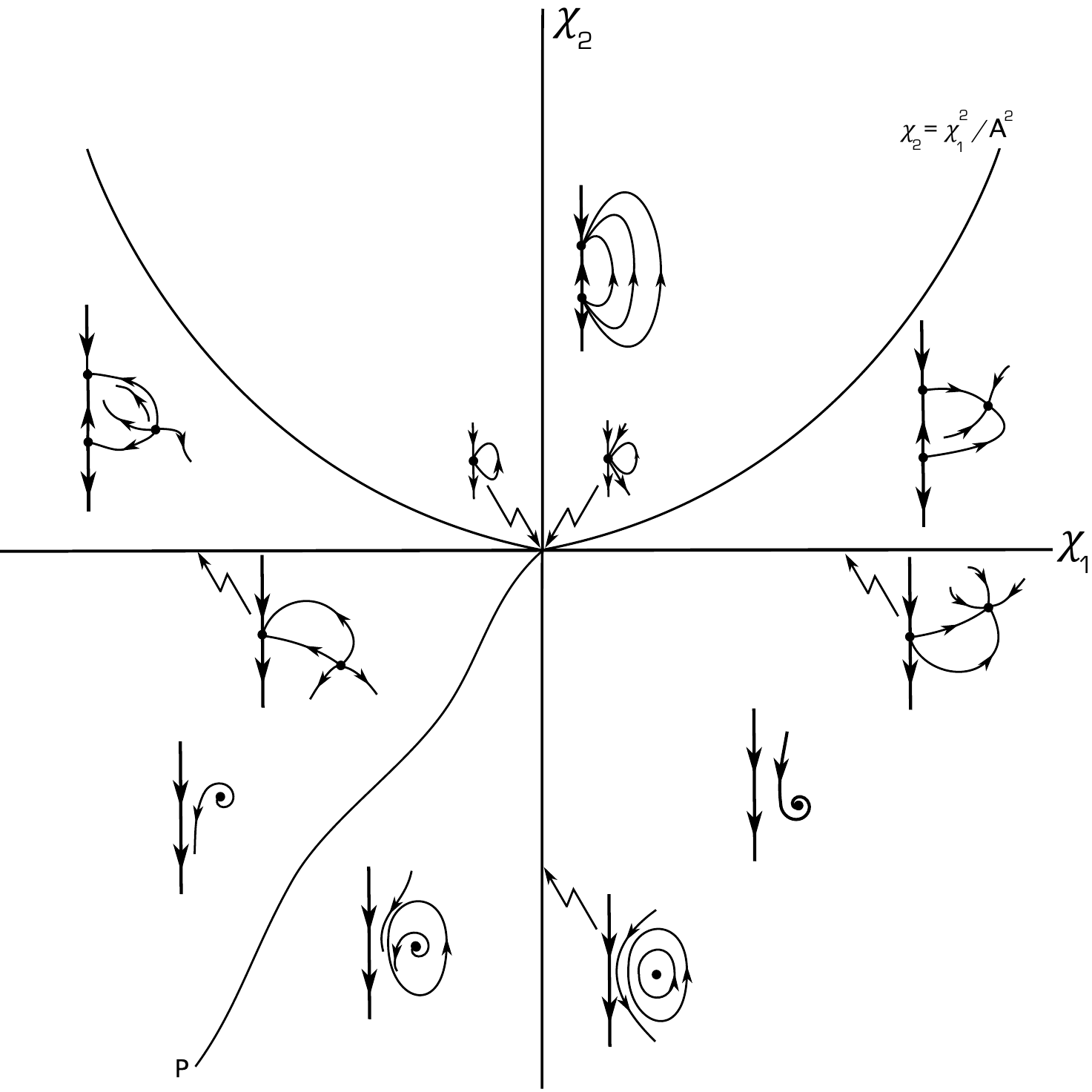}
\caption{Theoretical unfolding diagram of the Zero-Hopf bifuraction in the case $A<0$, $B=1$. Only the portion in the upper half-plane is attainable in (\ref{simsys2}).}
\label{fig1}
\end{center}
\end{figure}

The theroretical versal unfolding (in polar coordinates) for the Zero-Hopf bifurcation in the case $A<0$, $B=1$ is 
given in 
(\ref{NormalForm3}),
and the corresponding unfolding diagram is
illustrated in Figure \ref{fig1} (see \cite{Kuznetsov}), where we recall that according to Remark \ref{rmk42} (b), only the upper half-plane portion of this diagram is attainable in the unfolding of (\ref{simsys2}).  We remind the reader that the link between the theoretical unfolding parameters $\chi_{1,2}$ and the unfolding parameter $\mu_{1,2}$ in (\ref{simsys2}) are given by (\ref{ABdefs}).

\begin{figure}[htbp]
\begin{center}
\includegraphics[width=4in]{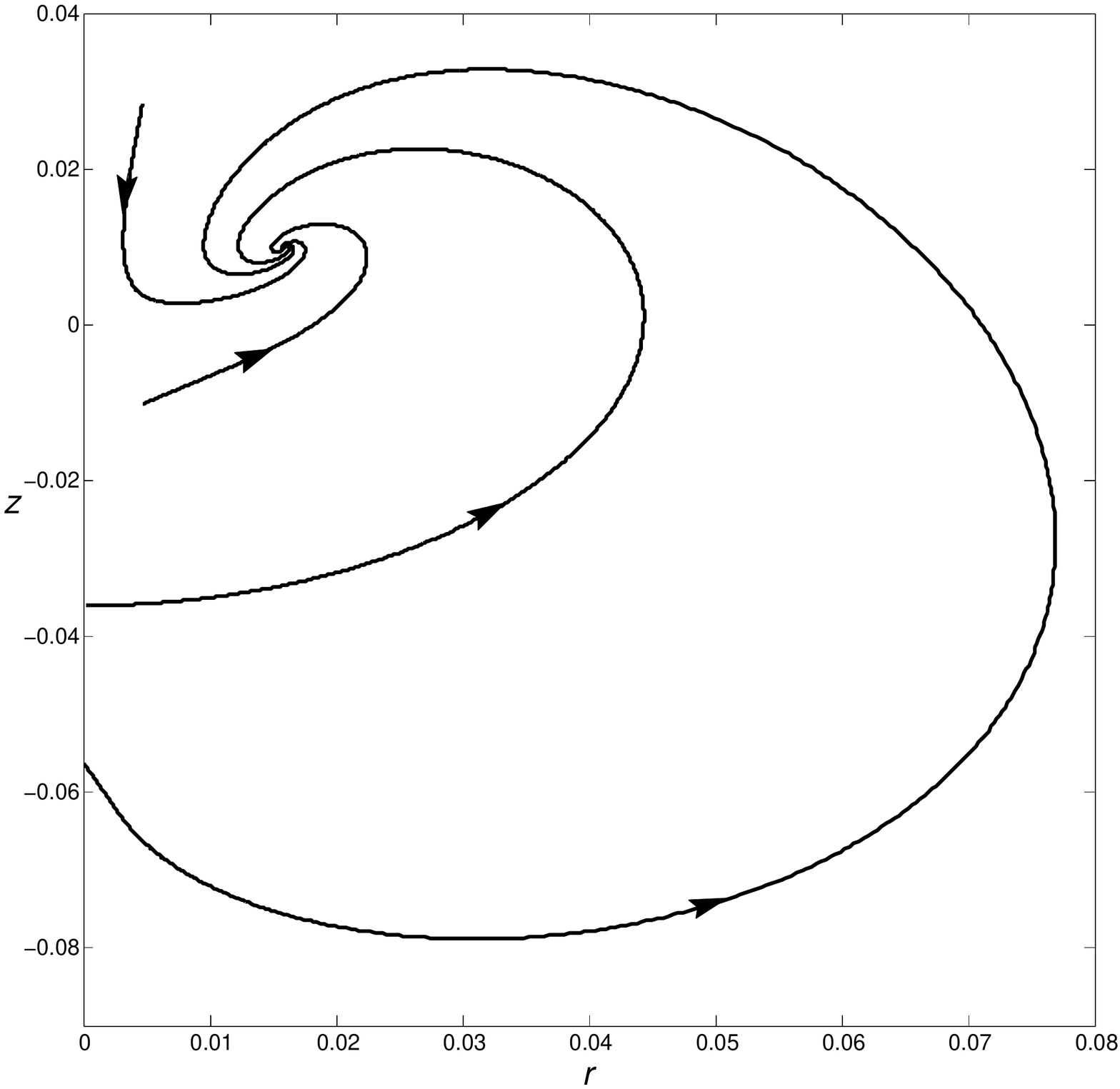}
\caption{Representation of numerical simulations of (\ref{simsys2}) with 4 different initial conditions, and parameter values (\ref{pm1}).  In this figure and in Figures \ref{fig4}, \ref{fig6}, \ref{fig8}, \ref{fig9} and \ref{fig10}, we have projected the numerical solution onto the eigenbasis of the center subspace using the bilinear form (\ref{BilinearForm}) and converted into polar coordinates $(r,z)$, which results in the illustrated trajectories.  As can be seen, trajectories tend towards an asymptotically stable non-trivial equilibrium point.}
\label{fig2}
\end{center}
\end{figure}

\begin{figure}[htbp]
\begin{center}
\includegraphics[width=4in]{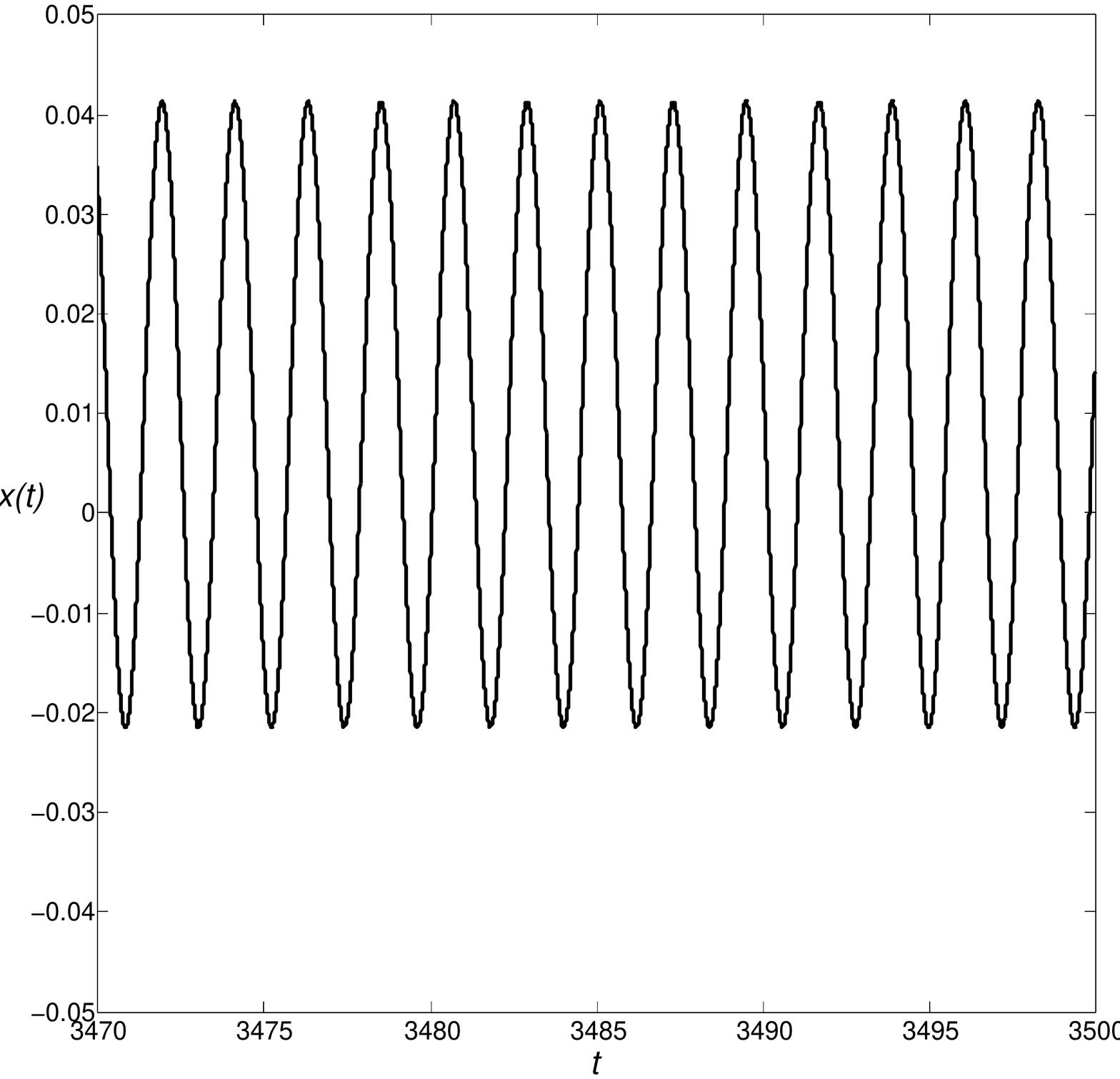}
\caption{The asymptotictically stable equilibrium point of figure \ref{fig2} corresponds to the time-periodic solution $x(t)$ of (\ref{simsys2}) depicted herein.}
\label{fig3}
\end{center}
\end{figure}

\begin{figure}[htbp]
\begin{center}
\includegraphics[width=4in]{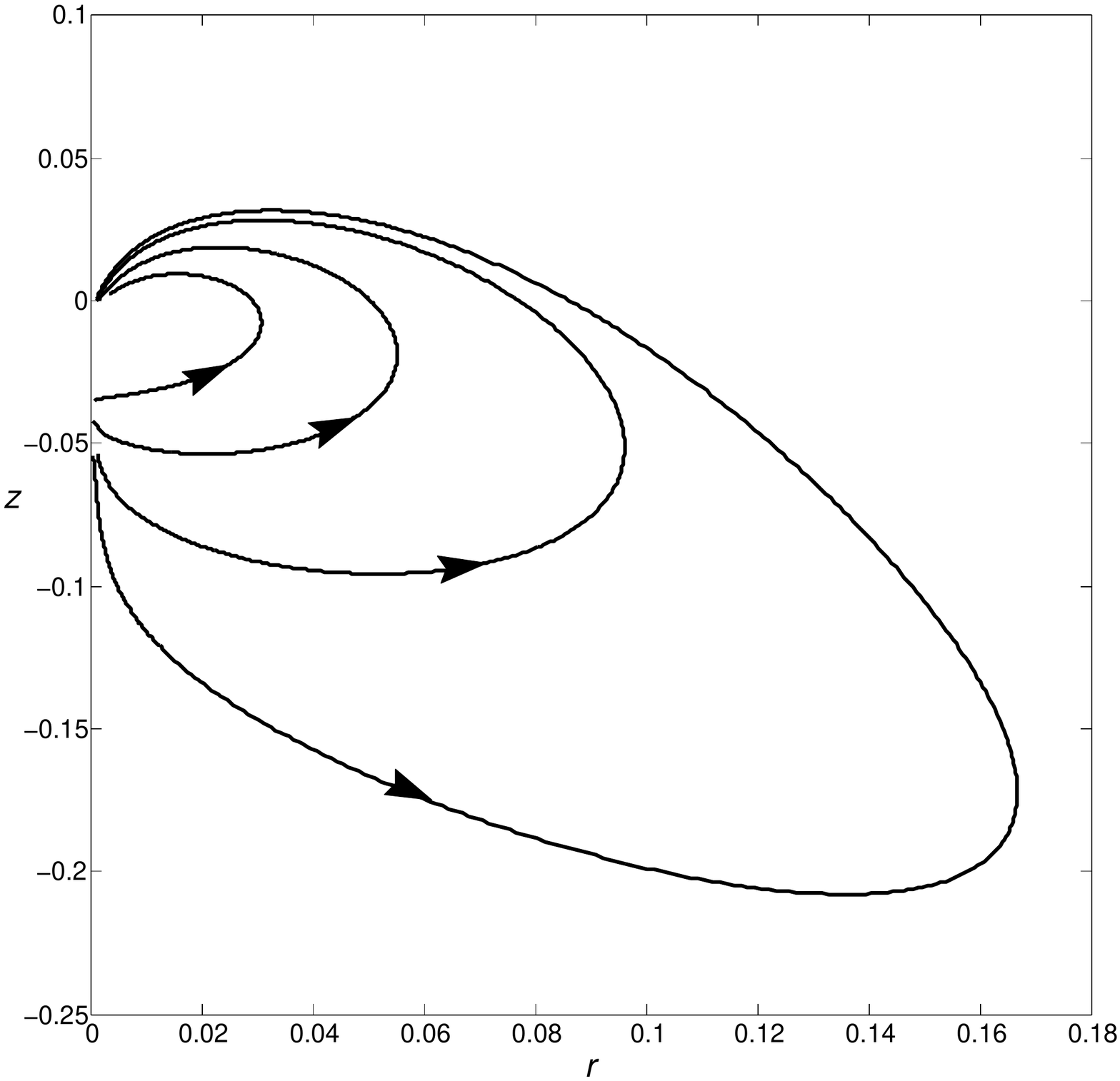}
\caption{Representation of numerical simulations of (\ref{simsys2}) with 4 different initial conditions, and parameter values (\ref{pm2}).   As can be seen, trajectories tend towards the asymptotically stable trivial equilibrium point.}
\label{fig4}
\end{center}
\end{figure}

In Figures \ref{fig2} and \ref{fig3}, we show the results of a numerical integration of (\ref{simsys2}) for parameter values 
\begin{equation}
\mu_1=-0.0018\,\,\,\,\,\,\,\,\,\,\,\,\,\,\,
\mu_2=0.0032.
\label{pm1}
\end{equation}
With respect to the theoretical unfolding diagram of Figure \ref{fig1}, this corresponds to a point $(\chi_1,\chi_2)$ in the first quadrant, below the parabola.
We note that the phase diagrams shown in Figures \ref{fig2}, \ref{fig4}, \ref{fig6}, \ref{fig8}, \ref{fig9} and \ref{fig10} were obtained by numerically integrating (\ref{simsys}) for different parameter values using {\em dde23} in {\em Matlab}, and then numerically computing the projection of the obtained solution $(x(t),\dot{x}(t))$ onto the center eigenspace $P$ of the Zero-Hopf bifurcation using the bilinear form (\ref{BilinearForm}), converting into polar coordinates $(r,z)$, and then using the {\em supsmu} package of Matlab to smooth out the trajectories.

In Figure \ref{fig4}, we show the results of a numerical integration of (\ref{simsys2}) for parameter values 
\begin{equation}
\mu_1=-0.0018\,\,\,\,\,\,\,\,\,\,\,\,\,\,\,
\mu_2=0.
\label{pm2}
\end{equation}
With respect to the theoretical unfolding diagram of Figure \ref{fig1}, this corresponds to a point $(\chi_1,\chi_2)$ above the parabola.

\subsection{Case II: $A<0$ and $B=-1$}

In this case, we simulate
\begin{equation}
\begin{array}{l}
\ddot{x}(t) + \varepsilon(x^2(t) - 1)\dot{x}(t) + x(t) = \\
\\
(1+\mu_1)x(t-\tau)+0.5\, \dot{x}(t-\tau)+0.2\, x(t-\tau)^2-0.1\, x(t-\tau)\dot{x}(t-\tau)+1\, \dot{x}(t-\tau)^2
\end{array}
\label{simsys3}
\end{equation}
with $\varepsilon=0.6$ and
where
\[
\tau=\tau_0+\mu_2.
\]

When $\mu_1=\mu_2=0$, system (\ref{simsys3}) 
undergoes a Zero-Hopf bifurcation with values:
\[
\omega_0\approx 1.375,\,\,\,\,\,\tau_0\approx 2.180
\]
as given in Proposition
\ref{EigenvalueTheorem} (v),
and
\[
A\approx -1.517,\,\,\,\,\,\,B=-1,
\]
as per (\ref{ABdefs}) and (\ref{aValues}).

\begin{figure}[htbp]
\begin{center}
\includegraphics[width=6in]{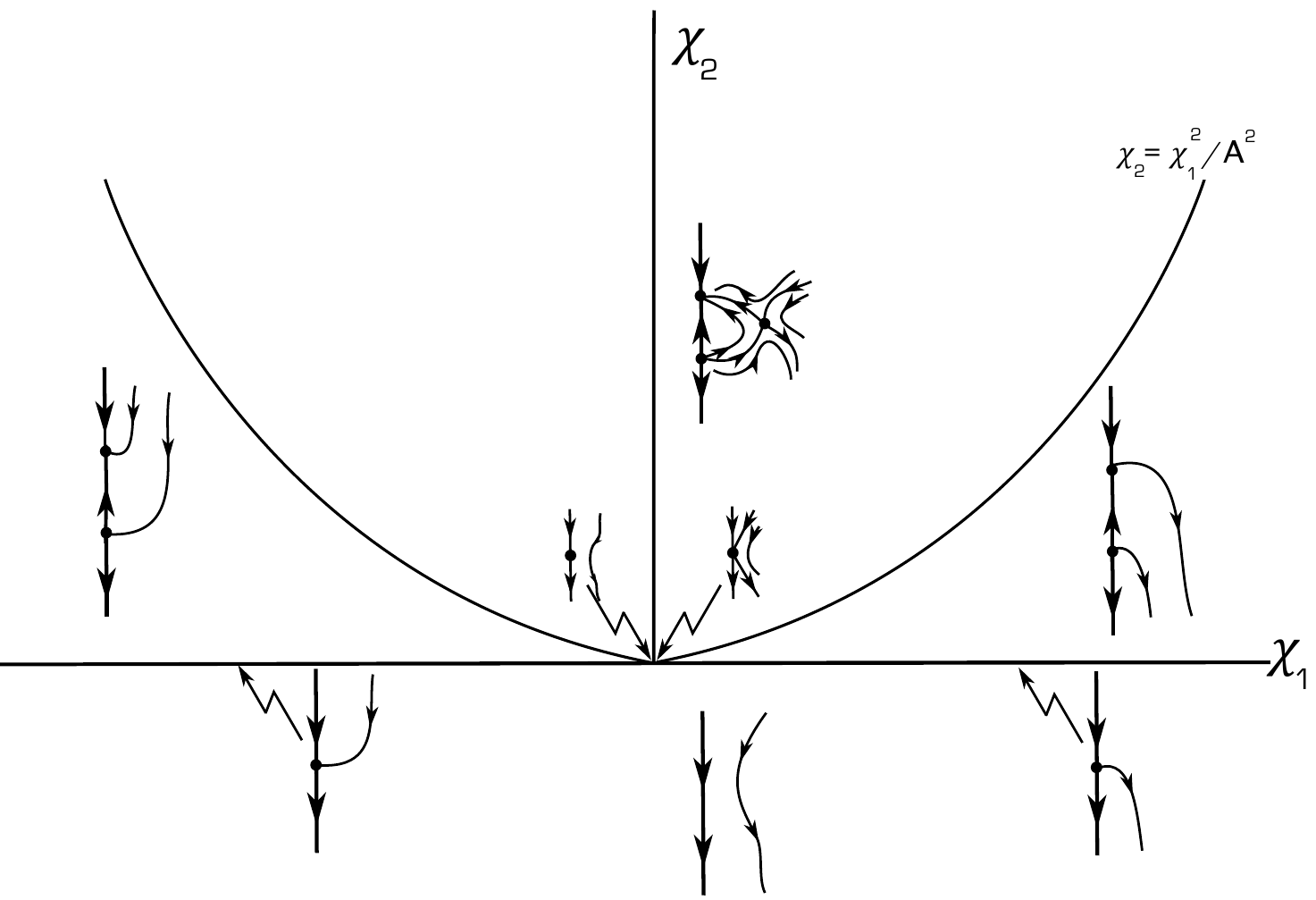}
\caption{Theoretical unfolding diagram of the Zero-Hopf bifuraction in the case $A<0$, $B=-1$. Only the portion in the upper half-plane is attainable in (\ref{simsys3}).}
\label{fig5}
\end{center}
\end{figure}

The theroretical versal unfolding (in polar coordinates) for the Zero-Hopf bifurcation in the case $A<0$, $B=-1$ is 
given in 
(\ref{NormalForm3}),
and the corresponding unfolding diagram is
illustrated in Figure \ref{fig5}.

\begin{figure}[htbp]
\begin{center}
\includegraphics[width=4in]{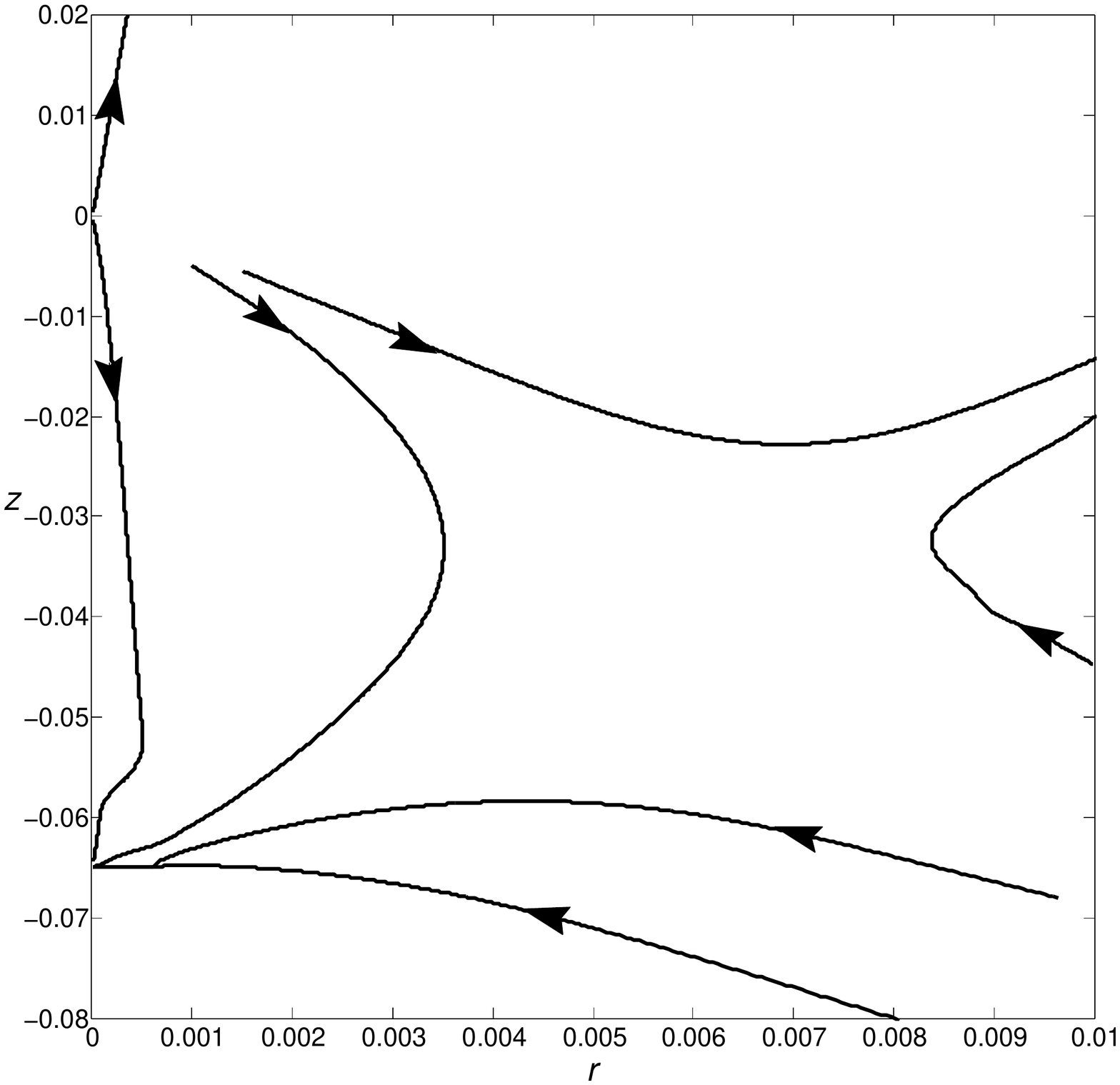}
\caption{Representation of numerical simulations of (\ref{simsys3}) with several different initial conditions, and parameter values (\ref{pm3}). 
We observe the behavior of the orbits in the vicinity of a saddle point, as well as the presence of both an asymptotically stable and an unstable trivial equilibrium point.}
\label{fig6}
\end{center}
\end{figure}

In Figure \ref{fig6}, we show the results of a numerical integration of (\ref{simsys3}) for parameter values 
\begin{equation}
\mu_1=0.00325\,\,\,\,\,\,\,\,\,\,\,\,\,\,\,
\mu_2=0.00192.
\label{pm3}
\end{equation}
With respect to the theoretical unfolding diagram of Figure \ref{fig5}, this corresponds to a point $(\chi_1,\chi_2)$ above the parabola.

\subsection{Case III: $A>0$ and $B=-1$}

In this final case, we simulate
\begin{equation}
\begin{array}{l}
\ddot{x}(t) + \varepsilon(x^2(t) - 1)\dot{x}(t) + x(t) = \\
\\
(1+\mu_1)x(t-\tau)-0.2\, \dot{x}(t-\tau)-0.2\, x(t-\tau)^2-0.2\, x(t-\tau)\dot{x}(t-\tau)-0.2\, \dot{x}(t-\tau)^2
\end{array}
\label{simsys4}
\end{equation}
with $\varepsilon=0.3$ and
where
\[
\tau=\tau_0+\mu_2.
\]

When $\mu_1=\mu_2=0$, system (\ref{simsys4}) 
undergoes a Zero-Hopf bifurcation with values:
\[
\omega_0\approx 1.396,\,\,\,\,\,\tau_0\approx 1.757
\]
as given in Proposition
\ref{EigenvalueTheorem} (v),
and
\[
A\approx 3.024,\,\,\,\,\,\,B=-1,
\]
as per (\ref{ABdefs}) and (\ref{aValues}).

\begin{figure}[htbp]
\begin{center}
\includegraphics[width=6in]{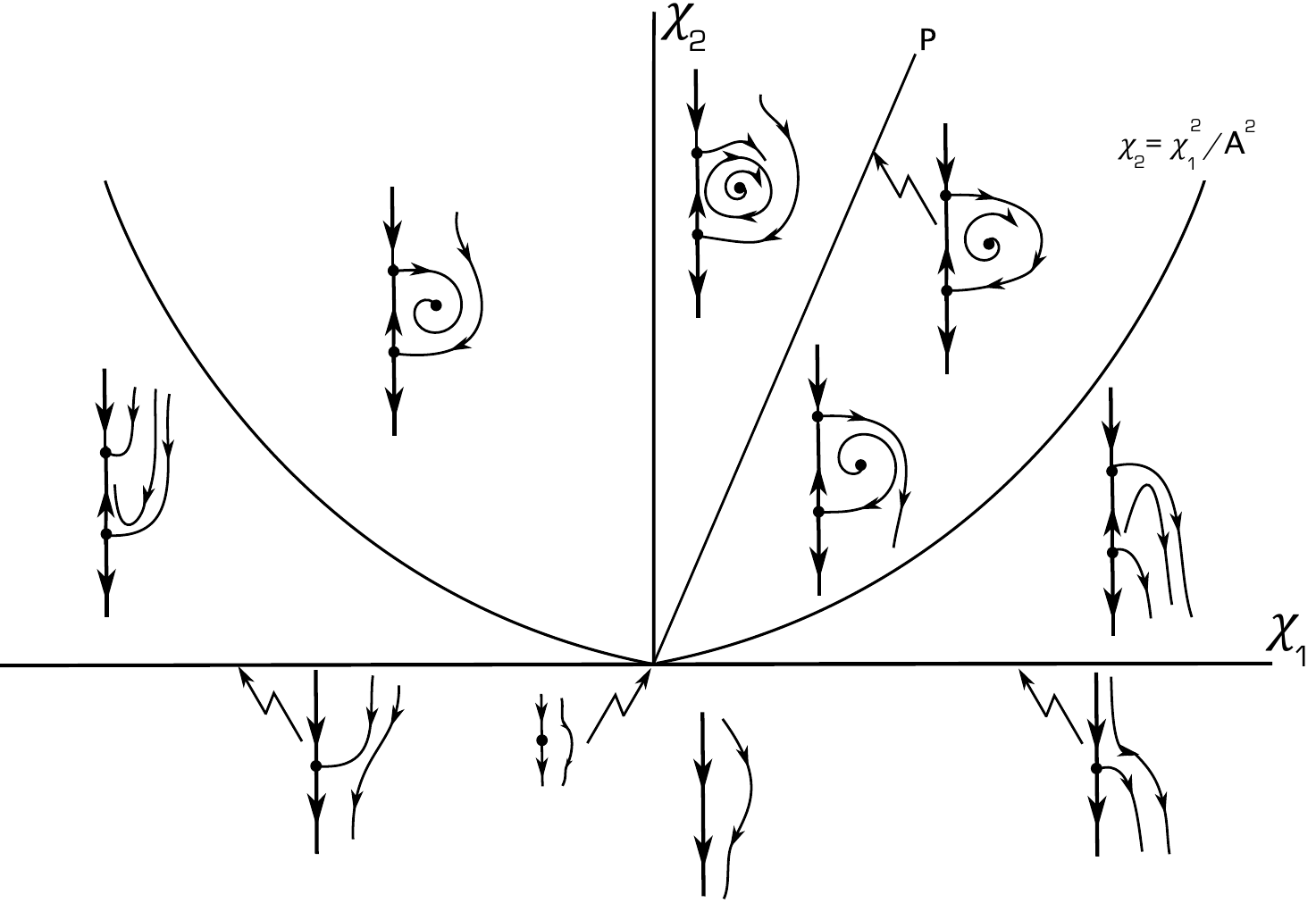}
\caption{Theoretical unfolding diagram of the Zero-Hopf bifuraction in the case $A>0$, $B=-1$. Only the portion in the upper half-plane is attainable in (\ref{simsys4}).}
\label{fig7}
\end{center}
\end{figure}

The theroretical versal unfolding diagram is
illustrated in Figure \ref{fig7}.

\begin{figure}[htbp]
\begin{center}
\includegraphics[width=4in]{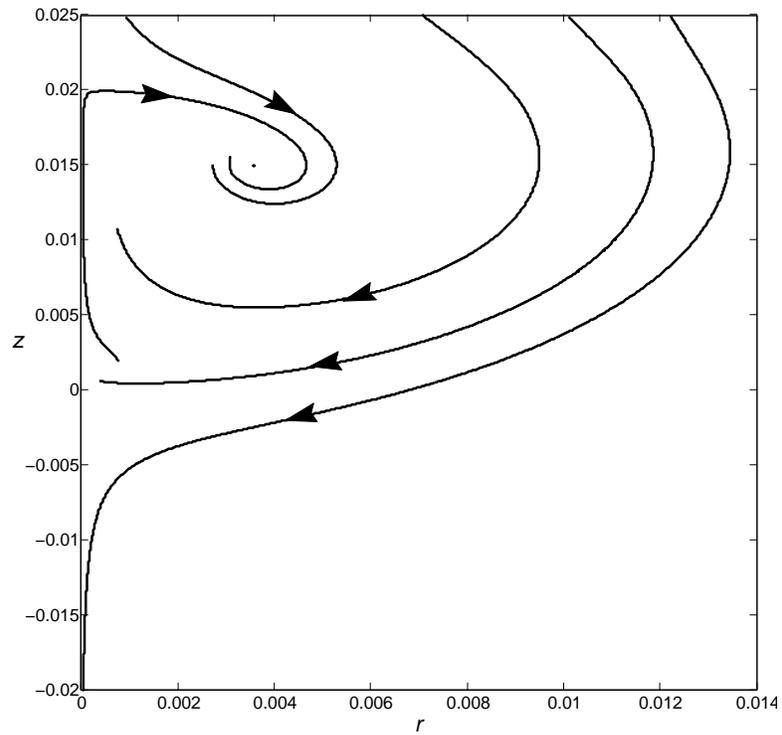}
\caption{Representation of numerical simulations of (\ref{simsys4}) with several different initial conditions, and parameter values (\ref{pm4}).   As can be seen, trajectories tend towards an asymptotically stable non-trivial equilibrium point, and there is a separatrix terminating at the saddle at the origin.}
\label{fig8}
\end{center}
\end{figure}

In Figure \ref{fig8}, we show the results of a numerical integration of (\ref{simsys4}) for parameter values 
\begin{equation}
\mu_1=0.001\,\,\,\,\,\,\,\,\,\,\,\,\,\,\,
\mu_2=-0.003.
\label{pm4}
\end{equation}
With respect to the theoretical unfolding diagram of Figure \ref{fig7}, this corresponds to a point $(\chi_1,\chi_2)$ in the second quadrant above the parabola.

\begin{figure}[htbp]
\begin{center}
\includegraphics[width=4in]{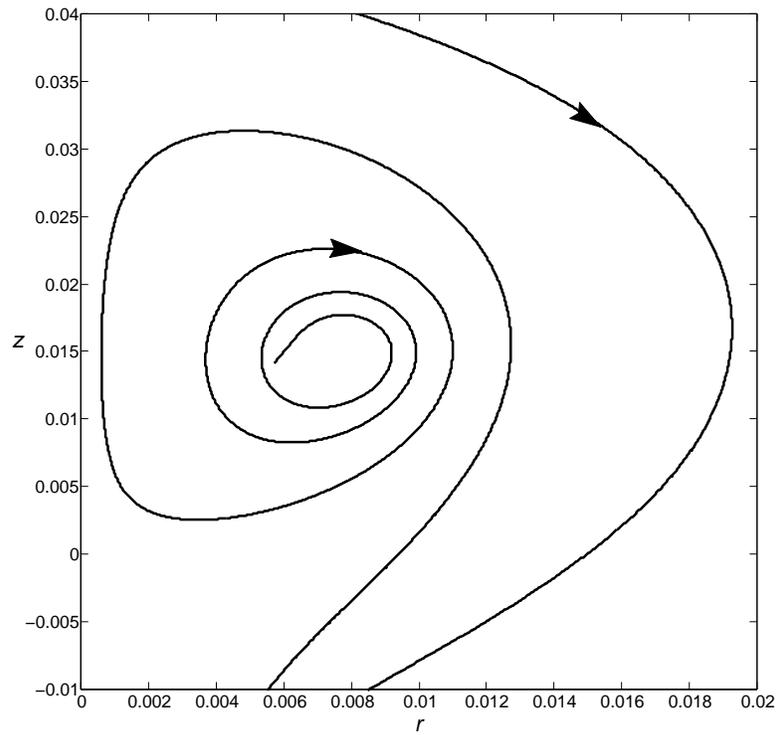}
\caption{Representation of numerical simulations of (\ref{simsys4}) with different initial conditions, and parameter values (\ref{pm5}).  The non-trivial equlilibrium point is a source.}
\label{fig9}
\end{center}
\end{figure}

\begin{figure}[htbp]
\begin{center}
\includegraphics[width=4in]{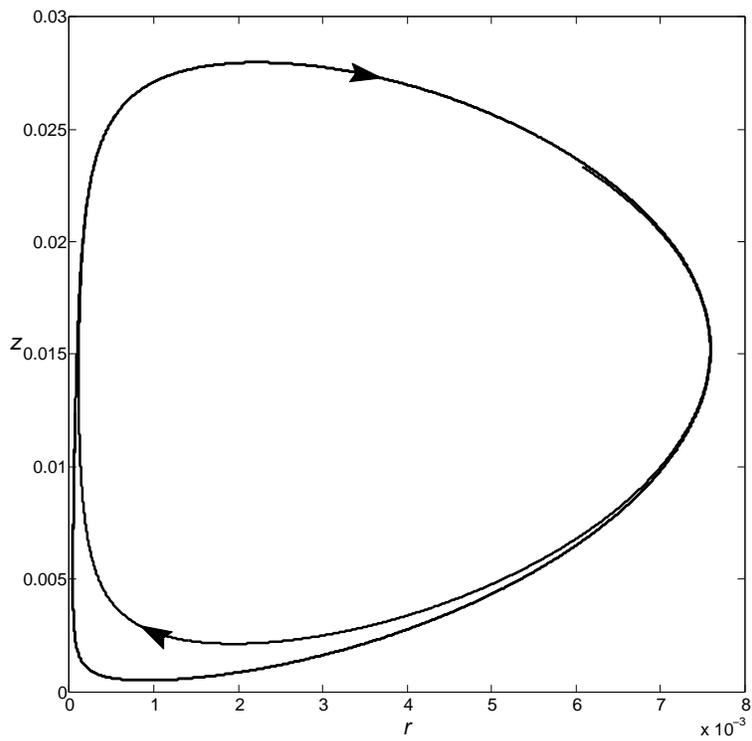}
\caption{Representation of a numerical simulation of (\ref{simsys4}) with parameter values (\ref{pm6}).    The trajectory appears to tend towards a limit cycle, but further numerical investigation would be required to establish this definitively.}
\label{fig10}
\end{center}
\end{figure}

In Figure \ref{fig9}, we show the results of a numerical integration of (\ref{simsys4}) for parameter values 
\begin{equation}
\mu_1=0.0019\,\,\,\,\,\,\,\,\,\,\,\,\,\,\,
\mu_2=-0.003.
\label{pm5}
\end{equation}
With respect to the theoretical unfolding diagram of Figure \ref{fig7}, this corresponds to a point $(\chi_1,\chi_2)$ in the first quadrant above the parabola and to the right of the heteroclinic bifurcation line $P$.  The non-trivial equilibrium point has gone from a sink in Figure \ref{fig8} to a source in Figure \ref{fig9}.   Of course, as can be seen in the theoretical bifurcation diagram in Figure \ref{fig7}, mediating this loss of stability is a Hopf bifurcation, followed by a heteroclinic bifurcation.  The limit cycle generated by the Hopf bifurcation grows in amplitude until it disappears in a heteroclinic connection between the two equilibrium points on the $r=0$ axis.  From (\ref{hetbifcurve}), we see that the Hopf bifurcation curve and the heteroclinic bifurcation curve are tangent at the origin.  Consequently, the region between the Hopf bifurcation curve and the heteroclinic bifurcation curve is very thin, and thus extremely hard to resolve numerically.  We have attempted to find parameter values $(\mu_1,\mu_2)$ in this region where (\ref{simsys4}) should have a stable limit cycle.  In Figure \ref{fig10}, we give the results of such an integration for
\begin{equation}
\mu_1=0.0015525\,\,\,\,\,\,\,\,\,\,\,\,\,\,\,
\mu_2=-0.003
\label{pm6}
\end{equation}
While the illustrated trajectory looks like it tends to a limit cycle, we can not conclusively say that it does.  Further detailed numerical analysis would be required to establish this definitively.

\appendix
\Section{Proof of Theorem \ref{MainThm}}

First, let $V_j(\mathbb{C}^3\ \times$ Ker $\pi)$ be the space of homogeneous polynomials of degree $j$ in the variables ($x,\mu$) with coefficients in $\mathbb{C}^3\ \times$ Ker $\pi$. Define $M_j: V_j(\mathbb{C}^3\ \times$ Ker $\pi) \to V_j(\mathbb{C}^3\ \times$ Ker $\pi)$ such that     
\begin{equation} \label{OperatorM}
	M_j(p,h) = (M^1_jp,M^2_jh)
\end{equation}
where
\begin{equation} \label{MDefinition}
	\begin{split}
	\begin{aligned}
		&M_j^1p(x,\mu) = D_xp(x,\mu)Jx - Jp(x,\mu) = i\omega \left( \begin{array}{cc}
 			x_1\dfrac{\partial p_1}{\partial x_1} - x_2\dfrac{\partial p_1}{\partial x_2} - p_1 \\ [0.4 cm]
 			x_1\dfrac{\partial p_2}{\partial x_1} - x_2\dfrac{\partial p_2}{\partial x_2} + p_2 \\ [0.4 cm]
			x_1\dfrac{\partial p_3}{\partial x_1} - x_2 \dfrac{\partial p_3}{\partial x_2}
 		\end{array}  \right), \\
		&M_j^2h(x,\mu) = D_xh(x,\mu)Jx - \mathcal{A}_{Q^1}h(x,\mu),
	\end{aligned}
	\end{split}
\end{equation}
with $p(x,\mu) \in V_j(\mathbb{C}^3)$ and $h(x,\mu)(\theta) \in V_j($Ker $\pi)$. We now have the decomposition $V_j(\mathbb{C}^3) = {\rm Im}(M^1_j) \otimes {\rm Ker}(M^1_j)$ for $j \geq 2$.

Introduce the nonlinear change of coordinates 
\begin{equation}
	\begin{split}
	\begin{aligned}
		&x = \hat{x} + U_2^1(\hat{x},\mu) \\
		&y = \hat{y} + U_2^2(\hat{x},\mu),
	\end{aligned}
	\end{split} 
\end{equation}
where  
\begin{equation} \label{ChangeOfCoordinates}
\begin{split}
\begin{aligned}
	U_2^1(x,\mu) &= (M_1^1)^{-1}Proj_{{\rm Im}(M_2^1)}f_2^1(x,0,\mu), \\
	U_2^2(x,\mu) &= (M_2^2)^{-1}f_2^2(x,0,\mu).
\end{aligned}
\end{split}
\end{equation}
Equation (\ref{ExpandedDecomposedSystem}) now becomes 
\begin{equation} 
	\begin{split}
	\begin{aligned}
		&\dot{x} = (I + D_xU_2^1(x,\mu))^{-1}\big[Jx + JU_2^1(x,\mu) + \sum_{j \geq 2}f_j^1(x+U_2^1(x,\mu),y+U_2^2(x,\mu))\big] \\
		&\frac{d}{dt}y = \mathcal{A}_{Q^1}y + \mathcal{A}_{Q^1}U_2^2(x,\mu) - D_xU_2^2(x,\mu)\dot{x} +\sum_{j \geq 2} f_j^2(x+U_2^1(x,\mu),y+U_2^2(x,\mu))
	\end{aligned}
	\end{split} 	
\end{equation} 
upon dropping the hats.

For $|x|$ small we have that 
\begin{equation} \label{InverseExpansion}
	(I + D_xU_2^1(x,\mu))^{-1} \approx I - D_xU_2^1(x,\mu) + (D_xU_2^1(x,\mu))^2,
\end{equation}
and using Taylor's Theorem we obtain
\[
	f_2^1(x+U_2^1(x,\mu),y+U_2^2(x,\mu)) = f_2^1(x,y) + D_xf^1_2(x,y)U_2^1(x,\mu) + D_yf^1_2(x,y)U_2^2(x,\mu) + h.o.t., 
\]
\[
	f_3^1(x+U_2^1(x,\mu),y+U_2^2(x,\mu))	= f_3^1(x,y) + h.o.t..
\]
Therefore we now have the normal form on the centre manifold given by 
\begin{equation} \label{NormalFormCentreManifold}
	\dot{x} = Jx + g^1_2(x,0,\mu) + g^1_3(x,0,\mu) + \mathcal{O}(|\mu|^2 + |x||\mu|^2 + |x|^2|\mu| + |x|^3),
\end{equation}  
where $g^1_2$ and $g^1_3$ are the following second and third order terms in ($x,\mu$): 
\begin{equation}
	\begin{split}
	\begin{aligned}
		g_2^1(x,0,\mu) = &Proj_{Ker(M_2^1)}f_2^1(x,0,\mu) \\
		g_3^1(x,0,\mu) = &Proj_{Ker(M_3^1)}\bigg[f_3^1(x,0,\mu) + D_xf_2^1(x,0,\mu)U_2^1(x,\mu) + D_yf(x,0,\mu)U_2^2(x,\mu)\\ &- D_xU_2^1(x,\mu)JU_2^1(x,\mu) - D_xU_2^1(x,\mu)f_2^1(x,0,\mu) + (D_xU_2^1(x,\mu))^2Jx\bigg].	
	\end{aligned}
	\end{split}
\end{equation}

Consider the basis $\{\mu^px^qe_k: k=1,2,3, p \in \mathbb{N}^2_0, q \in \mathbb{N}^3_0, |p| + |q| = j\}$ of $V_j(\mathbb{C}^3)$, where $e_1=(1,0,0)^T$, $e_2 = (0,1,0)^T$ and $e_3 = (0,0,1)^T$. Then for $j=2$, upon finding the images of each basis element under $M_2^1$, we find that Ker($M_2^1$) is spanned by 
\[
	\mu_1x_1e_1, \mu_2x_1e_1, x_1x_3e_1,
\]
\[
	\mu_1x_2e_2, \mu_2x_2e_2, x_2x_3e_2,
\]
\[
	x_1x_2e_3, \mu_1x_3e_3, \mu_2x_3e_3, \mu_1^2e_3, \mu_2^2e_3, \mu_1\mu_2e_3, x_3^2e_3.	
\]
Similarly, we find that Ker($M_3^1$) is spanned by
\[
	\mu_1\mu_2x_1e_1, \mu_1^2x_1e_1, \mu_2^2x_1e_1, \mu_1x_1x_3e_1, \mu_2x_1x_3e_1, x_1^2x_2e_1, x_1x_3^2e_1,
\]
\[	
	\mu_1\mu_2x_2e_2, \mu_1^2x_2e_2, \mu_2^2x_2e_2, \mu_1x_2x_3e_2, \mu_2x_2x_3e_2, x_1x_2^2e_2, x_2x_3^2e_2, 
\]
\[	
	\mu_1x_1x_2e_3, \mu_2x_1x_2e_3, \mu_1x_3^2e_3, \mu_2x_3^2e_3, \mu_1\mu_2x_3e_3, x_1x_2x_3e_3, x_3^3e_3.
\]
To see the images of each of these individual basis elements under the operator $M_j^1$ see $\cite{WuWang}$. 

First we compute $g_2^1(x,0,\mu)$ in (\ref{NormalFormCentreManifold}).  We have that
\begin{equation} \label{g2}
\begin{split}
\begin{aligned}
	g_2^1(x,0,\mu) &= Proj_{{\rm Ker}(M_2^1)}f_2^1(x,0,\mu) + \mathcal{O}(|\mu|^2) \\
	&= \begin{pmatrix}
		(a_{11}\mu_1 + a_{12}\mu_2)x_1 + a_{13}x_1x_3 \\
		(\bar{a}_{11}\mu_1 + \bar{a}_{12}\mu_2)x_2 + \bar{a}_{13}x_2x_3 \\
		(a_{21}\mu_1 + a_{22}\mu_2)x_3 + a_{23}x_1x_2 + a_{24}x_3^2
	\end{pmatrix}
	+ \mathcal{O}(|\mu|^2).
\end{aligned}
\end{split}
\end{equation}
Now expanding out the terms in ($\ref{TaylorTerms1}$) we get
\begin{equation} \label{TaylorTerms2}
 \begin{split}
 \begin{aligned}
 	F^1_2 = &\mu_2(i\omega_0x_1 - i\omega_0x_2 + y_2(0)), \ \ \ \ \ F^1_3 = 0, \\
	F^2_2 = &-\mu_2(x_1 + x_2 + x_3 + y_1(0)) + \varepsilon\mu_2(i\omega_0x_1 - i\omega_0x_2 + y_2(0)) + a\mu_2(i\omega_0e^{-i\omega_0\tau_0}x_1\\ &- i\omega_0e^{i\omega_0\tau_0}x_2 + y_2(-1)) + \tau_0\mu_1(e^{-i\omega_0\tau_0}x_1 + e^{i\omega_0\tau_0}x_2 + x_3 + y_1(-1)) +\mu_2(e^{-i\omega_0\tau_0}x_1\\ &+ e^{i\omega_0\tau_0}x_2 + x_3 + y_1(-1)) + \frac{1}{2}\tau_0g_{u_1u_1}(0,0)(e^{-i\omega_0\tau_0}x_1 + e^{i\omega_0\tau_0}x_2 + x_3 + y_1(-1))^2\\ &+ \tau_0g_{u_1u_2}(0,0)(e^{-i\omega_0\tau_0}x_1 + e^{i\omega_0\tau_0}x_2 + x_3 + y_1(-1))(i\omega_0e^{-i\omega_0\tau_0}x_1 - i\omega_0e^{i\omega_0\tau_0}x_2\\ &+ y_2(-1)) + \frac{1}{2}\tau_0g_{u_2u_2}(0,0)(i\omega_0e^{-i\omega_0\tau_0}x_1 - i\omega_0e^{i\omega_0\tau_0}x_2 + y_2(-1))^2, \\
	F_3^2 = &-\varepsilon\tau_0(x_1 + x_2 + x_3 + y_1(0))^2(i\omega_0x_1 - i\omega_0x_2 + y_2(0)) + \frac{1}{2}g_{u_1u_1}(0,0)\mu_2(e^{-i\omega_0\tau_0}x_1\\ &+ e^{i\omega_0\tau_0}x_2 + x_3 + y_1(-1))^2 + g_{u_1u_2}(0,0)\mu_2(e^{-i\omega_0\tau_0}x_1 + e^{i\omega_0\tau_0}x_2 + x_3 + y_1(-1))\\ &(i\omega_0e^{-i\omega_0\tau_0}x_1 - i\omega_0e^{i\omega_0\tau_0}x_2 + y_2(-1))+ \frac{1}{2}g_{u_2u_2}(0,0)\mu_2(i\omega_0e^{-i\omega_0\tau_0}x_1 - i\omega_0e^{i\omega_0\tau_0}x_2\\ &+ y_2(-1)) + \frac{1}{6}\tau_0g_{u_1u_1u_1}(0,0)(e^{-i\omega_0\tau_0}x_1+ e^{i\omega_0\tau_0}x_2 + x_3 + y_1(-1))^3 + \frac{1}{2}\tau_0g_{u_1u_1u_2}(0,0)\\ &(e^{-i\omega_0\tau_0}x_1+ e^{i\omega_0\tau_0}x_2 + x_3 + y_1(-1))^2(i\omega_0e^{-i\omega_0\tau_0}x_1 - i\omega_0e^{i\omega_0\tau_0}x_2 + y_2(-1))\\ &+ \frac{1}{2}\tau_0g_{u_2u_2u_1}(0,0)(e^{-i\omega_0\tau_0}x_1+ e^{i\omega_0\tau_0}x_2 + x_3 + y_1(-1))(i\omega_0e^{-i\omega_0\tau_0}x_1 - i\omega_0e^{i\omega_0\tau_0}x_2\\ &+ y_2(-1))^2 + \frac{1}{6}\tau_0g_{u_2u_2u_2}(0,0)(i\omega_0e^{-i\omega_0\tau_0}x_1 - i\omega_0e^{i\omega_0\tau_0}x_2 + y_2(-1))^3.
 \end{aligned}
 \end{split}
\end{equation}
Therefore, recalling that the characteristic equation gives us $(1+ia\omega_0)e^{-i\omega_0\tau_0} = -\omega_0^2 - i\varepsilon\omega_0 + 1$, we compute the coefficients in ($\ref{g2}$) to be
\begin{equation} \label{aValues} 
	\begin{split}
	\begin{aligned}
		&a_{11} = \tau_0\bar{D}\bar{\sigma}e^{-i\omega_0\tau_0}, \ \ \ a_{12} = \bar{D}(i\omega_0 - \bar{\sigma}\omega_0^2), \\
		&a_{13} = \tau_0\bar{D}\bar{\sigma}(g_{u_1u_1}(0,0)e^{-i\omega_0\tau_0} + i\omega_0g_{u_1u_2}(0,0)e^{-i\omega_0\tau_0}), \ \ \ a_{21} = \frac{\tau_0}{\tau_0 - \varepsilon - a}, \\
		&a_{22} = 0, \ \ \ a_{23} = a_{21}(g_{u_1u_1}(0,0) + \omega_0^2g_{u_2u_2}(0,0)), \ \ \ a_{24} = \frac{1}{2}a_{21}g_{u_1u_1}(0,0).
	\end{aligned}
	\end{split}	
\end{equation}

Now to compute $g_3^1(x,0,0)$, recall that 
\begin{equation} \label{g3}
\begin{split}
\begin{aligned} 
	g_3^1(x,0,\mu) &= Proj_{Ker(M_3^1)}\bigg[f_3^1(x,0,\mu) + D_xf_2^1(x,0,\mu)U_2^1(x,\mu) + D_yf(x,0,\mu)U_2^2(x,\mu)\\ &- D_xU_2^1(x,\mu)JU_2^1(x,\mu) - D_xU_2^1(x,\mu)f_2^1(x,0,\mu) + (D_xU_2^1(x,\mu))^2Jx\bigg] \\
	&= Proj_{{\rm Ker}(M_3^1)}f_3^1(x,0,0) + Proj_{{\rm Ker}(M_3^1)}D_xf_2^1(x,0,\mu)U_2^1(x,0)\\ &+ Proj_{{\rm Ker}(M_3^1)}D_yf(x,0,\mu)U_2^2(x,0) - Proj_{{\rm Ker}(M_3^1)}D_xU_2^1(x,\mu)JU_2^1(x,\mu)\\ &- Proj_{{\rm Ker}(M_3^1)}D_xU_2^1(x,\mu)f_2^1(x,0,0) + Proj_{{\rm Ker}(M_3^1)}(D_xU_2^1(x,0))^2Jx\\ &+ \mathcal{O}(|x||\mu|^2 + |x|^2|\mu|).
\end{aligned}
\end{split}	
\end{equation}
($i$) First we compute $Proj_{{\rm Ker}(M_3^1)}f_3^1(x,0,0)$. Since
\begin{equation}
	f_3^1(x,0,0) = \tau_0 \begin{pmatrix}
	\bar{D}\bar{\sigma}H_1(x_1,x_2,x_3) \\
	D\sigma H_1(x_1,x_2,x_3) \\
	-D_1 H_1(x_1,x_2,x_3)
	\end{pmatrix},	
\end{equation} 
where
\begin{equation}
	\begin{split}
 	\begin{aligned}
 	H_1(x_1,x_2,x_3) = &-i\varepsilon\omega_0(x_1 - x_2)(x_1 + x_2 + x_3)^2 + \frac{1}{6}g_{u_1u_1u_1}(0,0)(e^{-i\omega_0\tau_0}x_1\\ &+ e^{i\omega_0\tau_0}x_2 + x_3)^3 + \frac{1}{2}	g_{u_1u_1u_2}(0,0)(e^{-i\omega_0\tau_0}x_1 + e^{i\omega_0\tau_0}x_2 + x_3)^2\\ &(i\omega_0e^{-i\omega_0\tau_0}x_1 - i\omega_0e^{i\omega_0\tau_0}x_2) + \frac{1}{2}			g_{u_2u_2u_1}(0,0)(e^{-i\omega_0\tau_0}x_1 + e^{i\omega_0\tau_0}x_2\\ &+ x_3)(i\omega_0e^{-i\omega_0\tau_0}x_1 - i\omega_0e^{i\omega_0\tau_0}x_2)^2 - \frac{i\omega_0}{6}(e^{-i	\omega_0\tau_0}x_1 - e^{i\omega_0\tau_0}),
 	\end{aligned}
 	\end{split}
\end{equation}
this gives us that 
\begin{equation}
Proj_{{\rm Ker}(M_2^1)}f_3^1(x,0,0) = \begin{pmatrix}
	b_{11}x_1^2x_2 + b_{12}x_1x_3^2 \\
	\bar{b}_{11}x_1x_2^2 + \bar{b}_{12}x_2x_3^2	\\
	b_{21}x_1x_2x_3 + b_{22}x_3^3	
\end{pmatrix},
\end{equation}
such that
\begin{equation}
\begin{split}
\begin{aligned} 
	b_{11} = &\tau_0\bar{D}\bar{\sigma}e^{-i\omega_0\tau_0}\bigg[-2i\varepsilon\omega_0e^{i\omega_0\tau_0} + \frac{1}{3}g_{u_1u_1u_1}(0,0) + \frac{1}{2}g_{u_1u_1u_2}(0,0)(2\omega_0^2 - i\omega_0)\\ &+ \frac{1}{2}g_{u_2u_2u_1}(0,0)(2\omega_0^2 + i\omega_0) - \frac{i\omega_0^3}{3}g_{u_2u_2u_2}(0,0) \bigg], \\
	b_{12} = &\tau_0\bar{D}\bar{\sigma}e^{-i\omega_0\tau_0}\bigg[-i\varepsilon\omega_0e^{i\omega_0\tau_0} + \frac{1}{3}g_{u_1u_1u_1}(0,0) + \frac{i\omega_0}{2}g_{u_1u_1u_2}(0,0)\bigg] , \\
	b_{21} = &-\tau_0D_1\bigg[g_{u_1u_1u_1}(0,0) + \frac{2\omega_0^2}{3}g_{u_2u_2u_1}(0,0)\bigg], \ \ \ \ \ b_{22} = \frac{-\tau_0}{6}D_1g_{u_1u_1u_1}(0,0).
\end{aligned}
\end{split}
\end{equation}
($ii$) To compute $Proj_{{\rm Ker}(M_3^1)}D_xf_2^1(x,0,\mu)U_2^1(x,0)$ we use 
\begin{equation} \label{U21}
 \begin{split}
 \begin{aligned}
 	&U_2^1(x,0) = U_2^1(x,\mu)|_{\mu=0} = (M_2^1)^{-1}Proj_{Im(M_2^1)}f_2^1(x,0,0) \\
	&= \frac{\tau_0}{i\omega_0} \begin{pmatrix}
		\bar{D}\bar{\sigma}\bigg[g_{u_1u_1}(0,0)(e^{-2i\omega_0\tau_0}x_1^2 - 2x_1x_2 - \frac{1}{3}e^{2i\omega_0\tau_0}x_2^2 - x_3^2 - e^{i\omega_0\tau_0}x_2x_3)\\ + g_{u_2u_2}(0,0)(-\omega_0^2e^{-2i\omega_0\tau_0}x_1^2 - 2\omega_0^2x_1x_2 + \frac{1}{3}\omega_0^2e^{2i\omega_0\tau_0}x_2^2)\\ + 2g_{u_1u_2}(0,0)(i\omega_0e^{-2i\omega_0\tau_0}x_1^2 + \frac{1}{3}i\omega_0e^{2i\omega_0\tau_0}x_2^2 + \frac{1}{2}i\omega_0e^{i\omega_0\tau_0}x_2x_3)\bigg] \\
		D\sigma\bigg[g_{u_1u_1}(0,0)(\frac{1}{3}e^{-2i\omega_0\tau_0}x_1^2 + 2x_1x_2 - e^{2i\omega_0\tau_0}x_2^2 + x_3^2 + e^{-i\omega_0\tau_0}x_1x_3)\\ + g_{u_2u_2}(0,0)(-\frac{1}{3}\omega_0^2e^{-2i\omega_0\tau_0}x_1^2 + 2\omega_0^2x_1x_2 + \omega_0^2e^{2i\omega_0\tau_0}x_2^2)\\ + 2g_{u_1u_2}(0,0)(\frac{1}{3}i\omega_0e^{-2i\omega_0\tau_0}x_1^2 - i\omega_0e^{2i\omega_0\tau_0}x_2^2 + \frac{1}{2}i\omega_0e^{-i\omega_0\tau_0}x_1x_3)\bigg] \\
		-D_1\bigg[g_{u_1u_1}(0,0)(\frac{1}{2}e^{-2i\omega_0\tau_0}x_1^2 - \frac{1}{2}e^{2i\omega_0\tau_0}x_2^2 + 2e^{-i\omega_0\tau_0}x_1x_3 - 2e^{i\omega_0\tau_0}x_2x_3)\\ + g_{u_2u_2}(0,0)(-\frac{1}{2}\omega_0^2e^{-2i\omega_0\tau_0}x_1^2 + \frac{1}{2}\omega_0^2e^{2i\omega_0\tau_0}x_2^2)\\ + 2g_{u_1u_2}(0,0)(\frac{1}{2}i\omega_0e^{-2i\omega_0\tau_0}x_1^2 + \frac{1}{2}i\omega_0e^{2i\omega_0\tau_0}x_2^2 + i\omega_0e^{-i\omega_0\tau_0}x_1x_3 + i\omega_0e^{i\omega_0\tau_0}x_2x_3)\bigg]
	\end{pmatrix}.
 \end{aligned}
 \end{split}
\end{equation}
This gives us
\begin{equation}
	Proj_{{\rm Ker}(M_3^1)}D_xf_2^1(x,0,0)U_2^1(x,0) = \begin{pmatrix}
		c_{11}x_1^2x_2 + c_{12}x_1x_3^2 \\
		\bar{c}_{11}x_1x_2^2 + \bar{c}_{12}x_2x_3^2 \\
		c_{21}x_1x_2x_3 + c_{22}x_3^3
	\end{pmatrix},	
\end{equation} 
where 
\begin{equation}
 \begin{split}
 \begin{aligned} 
 	c_{11} &= \frac{-\tau_0^2\bar{D}\bar{\sigma}}{3i\omega_0}\bigg[6\bar{D}\bar{\sigma}e^{-2i\omega_0\tau_0}\bigg((g_{u_1u_1}(0,0))^2 - \omega_0^4(g_{u_2u_2}(0,0))^2 + 2i\omega_0^3g_{u_1u_2}(0,0)g_{u_2u_2}(0,0)\\ &+ 2i\omega_0g_{u_1u_1}(0,0)g_{u_1u_2}(0,0)\bigg) + 2D\sigma\bigg(-7(g_{u_1u_1}(0,0))^2 - 10\omega_0^2g_{u_1u_1}(0,0)g_{u_2u_2}(0,0)\\ &- 4\omega_0^2(g_{u_1u_2}(0,0))^2 -7\omega_0^2(g_{u_2u_2}(0,0))^2\bigg) + 3D_1e^{-i\omega_0\tau_0}\bigg((g_{u_1u_1}(0,0))^2\\ &+ i\omega_0g_{u_1u_1}(0,0)g_{u_1u_2}(0,0) + i\omega_0^3g_{u_1u_2}(0,0)g_{u_2u_2}(0,0) - \omega_0^2g_{u_1u_1}(0,0)g_{u_2u_2}(0,0)\\ &+ 2\omega_0^2(g_{u_1u_2}(0,0))^2\bigg)\bigg], \\
	c_{12} &= \frac{-2\tau_0\bar{D}\bar{\sigma}}{i\omega_0}\bigg[\bar{D}\bar{\sigma}e^{-2i\omega_0\tau_0}\bigg((g_{u_1u_1}(0,0))^2 + 2i\omega_0g_{u_1u_1}(0,0)g_{u_1u_2}(0,0)\\ &- \omega_0^2g_{u_1u_1}(0,0)g_{u_2u_2}(0,0)\bigg) + D\sigma\bigg(-2(g_{u_1u_1}(0,0))^2 - \omega_0^2g_{u_1u_1}(0,0)g_{u_2u_2}(0,0)\\ &- \omega_0^2(g_{u_1u_2}(0,0))^2\bigg) + D_1e^{-i\omega_0\tau_0}\bigg((g_{u_1u_1}(0,0))^2 + 2i\omega_0g_{u_1u_1}(0,0)g_{u_1u_2}(0,0)\bigg)\bigg], \\
	c_{21} &= \frac{2\tau_0^2D_1}{i\omega_0}\bigg[\bar{D}\bar{\sigma}e^{-i\omega_0\tau_0}\bigg(3(g_{u_1u_1}(0,0))^2 + 3i\omega_0g_{u_1u_1}(0,0)g_{u_1u_2}(0,0)\\ &+ 2\omega_0^2(g_{u_1u_2}(0,0))^2 + \omega_0^2g_{u_1u_1}(0,0)g_{u_2u_2}(0,0) + 3i\omega_0^3g_{u_1u_2}(0,0)g_{u_2u_2}(0,0)\bigg)\\ &+ D\sigma e^{i\omega_0\tau_0}\bigg(-3(g_{u_1u_1}(0,0))^2 + 3i\omega_0g_{u_1u_1}(0,0)g_{u_1u_2}(0,0) - 2\omega_0^2(g_{u_1u_2}(0,0))^2\\ &- \omega_0^2g_{u_1u_1}(0,0)g_{u_2u_2}(0,0) + 3i\omega_0^3g_{u_1u_2}(0,0)g_{u_2u_2}(0,0)\bigg)\bigg], \\
	c_{22} &= \frac{2\tau_0D_1g_{u_1u_1}(0,0)}{i\omega_0}\bigg[\bar{D}\bar{\sigma}e^{-i\omega_0\tau_0}\bigg(g_{u_1u_1}(0,0) + i\omega_0g_{u_1u_2}(0,0)\bigg)\\ &+ D\sigma e^{i\omega_0\tau_0}\bigg(-g_{u_1u_1}(0,0) + i\omega_0g_{u_1u_2}(0,0)\bigg)\bigg].
\end{aligned}
\end{split}
\end{equation}
($iii$) To compute $Proj_{Ker(M_3^1)}D_yf(x,0,\mu)U_2^2(x,0)$ we define $h = h(x)(\theta) = U_2^2$, and write
\begin{equation}
	h(\theta) = \begin{pmatrix}
		h^{(1)}(\theta) \\
		h^{(2)}(\theta) 
	\end{pmatrix} = h_{200}x_1^2 + h_{020}x_2^2 + h_{002}x_3^2 + h_{110}x_1x_2 + h_{101}x_1x_3 + h_{011}x_2x_3,
\end{equation}
where $h_{200}, h_{020}, h_{002}, h_{110}, h_{101}, h_{011} \in Q^1$. Then we solve for the coefficients of $h$ using the fact that $(M_2^2h)(x) = f_2^2(x,0,0)$, or equivalently,
\begin{equation}
	D_xhJx - \mathcal{A}_{Q^1}(h) = (I - \pi)X_0F_2(\Phi x,0).	
\end{equation}   
Applying the definition of $\mathcal{A}$ and $\pi$ we obtain the following ordinary differential equation,
\begin{equation}
\begin{split}
 \begin{aligned}
 	&\dot{h} - D_xhJx = \Phi(\theta)\Psi(0)F_2(\Phi x,0), \\
	&\dot{h}(0) - Lh = F_2(\Phi x,0),
\end{aligned}
\end{split}	
\end{equation}
where $\dot{h}$ is the derivative of $h$ with respect to $\theta$. 
If we denote 
\begin{equation}
	F_2(\Phi x,0) = A_{200}x_1^2 + A_{020}x_2^2 + A_{002}x_3^2 + A_{110}x_1x_2 + A_{101}x_1x_3 + A_{011}x_2x_3,
\end{equation}
where $A_{ijk} \in \mathbb{C}^2$, $0 \leq i,j,k \leq 2, i + j + k = 2$, we are able to compare the coefficients of each monomial and obtain a differential equation for each individual coefficient of $h$. An inspection of the coefficients in $F_2(\Phi x,0)$ reveals that $\bar{h}_{020} = h_{200}$ and $\bar{h}_{011} = h_{101}$, and therefore we must only solve the following ordinary differential equations,   
\begin{equation} \label{h200} 
    	\begin{array}{ll}
       		\dot{h}_{200} - 2i\omega_0\tau_0h_{200} = \Phi(\theta)\Psi(0)A_{200}, \\ 
		\dot{h}_{200}(0) - L(h_{200}) = A_{200},
	\end{array}
\end{equation} 
\begin{equation} \label{h101}
    	\begin{array}{ll}
       		\dot{h}_{101} - i\omega_0\tau_0h_{101} = \Phi(\theta)\Psi(0)A_{101}, \\ 
		\dot{h}_{101}(0) - L(h_{101}) = A_{101},
	\end{array}
\end{equation} 
\begin{equation}  \label{h110}
    	\begin{array}{ll}
       		\dot{h}_{110} = \Phi(\theta)\Psi(0)A_{110}, \\ 
		\dot{h}_{110}(0) - L(h_{110}) = A_{110},
	\end{array}
\end{equation}
\begin{equation} \label{h002}
    	\begin{array}{ll}
       		\dot{h}_{002} = \Phi(\theta)\Psi(0)A_{002}, \\ 
		\dot{h}_{002}(0) - L(h_{002}) = A_{002}.
	\end{array}
\end{equation}
Solving these linear differential equations is standard \cite{WuWang}, and the details are omitted.
  
Since
\begin{equation}
	F_2(u_t,0) = \begin{pmatrix}
	0 \\
	\tau_0g_{u_1u_1}(0,0)u^2_1(-1) + \tau_0g_{u_2u_2}(0,0)u^2_2(-1) + 2\tau_0g_{u_1u_2}(0,0)u_1(-1)u_2(-1) 
	\end{pmatrix},
\end{equation}
we have
\begin{equation}
	D_yf_2^1|_{y=0,\mu=0}(h) = \frac{\tau_0}{2} \begin{pmatrix}
		\psi_{12}H_2(x_1,x_2,x_3)(h) \\
		\psi_{22}H_2(x_1,x_2,x_3)(h) \\
		\psi_{32}H_2(x_1,x_2,x_3)(h)
	\end{pmatrix},	
\end{equation}
where
\begin{equation}
 	\begin{split}
	\begin{aligned}
 		H_2(x_1,x_2,x_3)(h) = &g_{u_1u_1}(0,0)(e^{-i\omega_0\tau_0}x_1 + e^{i\omega_0\tau_0}x_2 + x_3)h^{(1)}(-1) +\\ &g_{u_2u_2}(0,0)(i\omega_0e^{-i\omega_0\tau_0}x_1- i\omega_0e^{i\omega_0\tau_0}x_2)h^{(2)}(-1) + g_{u_1u_2}(0,0)(e^{-i\omega_0\tau_0}x_1\\ &+ e^{i\omega_0\tau_0}x_2 + x_3)h^{(2)}(-1)+ g_{u_1u_2}(0,0)(i\omega_0e^{-i\omega_0\tau_0}x_1\\ &- i\omega_0e^{i\omega_0\tau_0}x_2)h^{(1)}(-1).
 	\end{aligned}
	\end{split}
\end{equation}
Therefore, 
\begin{equation} 
	Proj_{{\rm Ker}(M_3^1)}(D_yf_2^1(x,y,\mu))U_2^2(x,\mu)|_{y=0,\mu=0} = \begin{pmatrix}
		d_{11}x_1^2x_2 + d_{12}x_1x_3^2 \\
		\bar{d}_{11}x_1x_2^2 + \bar{d}_{12}x_2x_3^2 \\
		d_{21}x_1x_2x_3 + d_{22}x_3^3 
	\end{pmatrix},
\end{equation}
where
\begin{equation}
	 \begin{split}
 \begin{aligned}
 	d_{11} = &\frac{\tau_0\bar{D}\bar{\sigma}}{2}\bigg[g_{u_1u_1}(0,0)(e^{i\omega_0\tau_0}h_{200}^{(1)}(-1) + e^{-i\omega_0\tau_0}h_{110}^{(1)}(-1)) + g_{u_1u_2}(0,0)(e^{i\omega_0\tau_0}h_{200}^{(2)}(-1)\\ &+ e^{-i\omega_0\tau_0}h_{110}^{(2)}(-1) + i\omega_0e^{-i\omega_0\tau_0}h_{110}^{(1)}(-1) - i\omega_0e^{i\omega_0\tau_0}h_{200}^{(1)}(-1))\\ &+ i\omega_0g_{u_2u_2}(0,0)(e^{-i\omega_0\tau_0}h_{110}^{(2)}(-1) - e^{i\omega_0\tau_0}h_{200}^{(2)}(-1))\bigg], \\
	d_{12} = &\frac{\tau_0\bar{D}\bar{\sigma}}{2}\bigg[g_{u_1u_1}(0,0)(e^{-i\omega_0\tau_0}h_{002}^{(1)}(-1) + h_{101}^{(1)}(-1)) + i\omega_0g_{u_2u_2}(0,0)e^{-i\omega_0\tau_0}h_{002}^{(2)}(-1)\\ &+ g_{u_1u_2}(0,0)(e^{-i\omega_0\tau_0}h_{002}^{(2)}(-1) + h_{101}^{(2)}(-1) + i\omega_0e^{-i\omega_0\tau_0}h_{002}^{(1)}(-1))\bigg], \\
	d_{21} = &\frac{-\tau_0D_1}{2}\bigg[g_{u_1u_1}(0,0)(e^{-i\omega_0\tau_0}h_{011}^{(1)}(-1) + e^{i\omega_0\tau_0}h_{101}^{(1)}(-1) + h_{110}^{(1)}(-1))\\ &+ i\omega_0g_{u_2u_2}(0,0)(e^{-i\omega_0\tau_0}h_{011}^{(2)}(-1) - e^{i\omega_0\tau_0}h_{101}^{(2)}(-1)) + g_{u_1u_2}(0,0)(e^{i\omega_0\tau_0}h_{011}^{(2)}(-1)\\ &+ e^{i\omega_0\tau_0}h_{101}^{(2)}(-1) + h_{110}^{(2)}(-1) + i\omega_0e^{-i\omega_0\tau_0}h_{011}^{(1)}(-1) - i\omega_0e^{i\omega_0\tau_0}h_{101}^{(1)}(-1))\bigg], \\
	d_{22} = &\frac{-\tau_0D_1}{2}\bigg[g_{u_1u_1}(0,0)h_{002}^{(1)}(-1) + g_{u_1u_2}(0,0)h_{002}^{(2)}(-1)\bigg].
 \end{aligned}
 \end{split}
\end{equation}

($iv$) Lastly, to compute $Proj_{Ker(M_3^1)}D_xU_2^1(x,\mu)JU_2^1(x,\mu)$, $\\Proj_{Ker(M_3^1)}D_xU_2^1(x,\mu)f_2^1(x,0,0)$ and $Proj_{Ker(M_3^1)}(D_xU_2^1(x,0))^2Jx$ we must only use the explicit form of $U_2^1(x,\mu)$ as given in $(\ref{U21})$.   

Then, 
\begin{equation}
	Proj_{{\rm Ker}(M_3^1)}D_xU_2^1(x,0)JU_2^1(x,0) = \begin{pmatrix}
		e_{11}x_1^2x_2 + e_{12}x_1x_3^2 \\
		\bar{e}_{11}x_1x_2^2 + \bar{e}_{12}x_2x_3^2 \\
		e_{21}x_1x_2x_3 + e_{22}x_3^3 
	\end{pmatrix},
\end{equation}
where
\begin{equation} 
	\begin{split}
\begin{aligned}
	e_{11} = &\frac{-2\tau_0^2\bar{D}\bar{\sigma}}{9i\omega_0}\bigg[\bar{D}\bar{\sigma}e^{-i\omega_0\tau_0}(27(g_{u_1u_1}(0,0))^2e^{-i\omega_0\tau_0} - 27\omega_0^4(g_{u_2u_2}(0,0))^2e^{-i\omega_0\tau_0}\\ &+ 18i\omega_0g_{u_1u_1}(0,0)g_{u_1u_2}(0,0)e^{-i\omega_0\tau_0} + 18i\omega_0^3g_{u_1u_2}(0,0)g_{u_2u_2}(0,0)e^{-i\omega_0\tau_0} +\\ &36i\omega_0g_{u_1u_1}(0,0)g_{u_1u_2}(0,0) + 36i\omega_0^3g_{u_1u_2}(0,0)g_{u_2u_2}(0,0)) + D\sigma(-19(g_{u_1u_1}(0,0))^2\\ &- 34\omega_0^2g_{u_1u_1}(0,0)g_{u_2u_2}(0,0) - 19\omega_0^4(g_{u_2u_2}(0,0))^2 - 4\omega_0^2(g_{u_1u_2}(0,0))^2)\bigg], \\
	e_{12} = &\frac{-\tau_0^2\bar{D}\bar{\sigma}}{i\omega_0}\bigg[\bar{D}\bar{\sigma}e^{-i\omega_0\tau_0}(2(g_{u_1u_1}(0,0))^2e^{-i\omega_0\tau_0} + 4i\omega_0g_{u_1u_1}(0,0)g_{u_1u_2}(0,0)\\ &- 2\omega_0^2g_{u_1u_1}(0,0)g_{u_2u_2}(0,0)e^{-i\omega_0\tau_0}) + D\sigma(-3(g_{u_1u_1}(0,0))^2\\ &- 2\omega_0^2g_{u_1u_1}(0,0)g_{u_2u_2}(0,0) - \omega_0^2(g_{u_1u_2}(0,0))^2)\bigg],\\
	e_{21} = &\frac{\tau_0^2D_1}{i\omega_0}\bigg[\bar{D}\bar{\sigma}e^{-i\omega_0\tau_0}(5(g_{u_1u_1}(0,0))^2 + 3\omega_0^2g_{u_1u_1}(0,0)g_{u_2u_2}(0,0)\\ &+ 5i\omega_0^3g_{u_1u_2}(0,0)g_{u_2u_2}(0,0) + 5i\omega_0g_{u_1u_1}(0,0)g_{u_1u_2}(0,0) + 2\omega_0^2(g_{u_1u_2}(0,0))^2)\\ &- D\sigma e^{i\omega_0\tau_0}(5(g_{u_1u_1}(0,0))^2 + 3\omega_0^2g_{u_1u_1}(0,0)g_{u_2u_2}(0,0) - 5i\omega_0^3g_{u_1u_2}(0,0)g_{u_2u_2}(0,0)\\ &- 5i\omega_0g_{u_1u_1}(0,0)g_{u_1u_2}(0,0) + 2\omega_0^2(g_{u_1u_2}(0,0))^2)\bigg],\\
	e_{22} = &\frac{2\tau_0^2D_1g_{u_1u_1}(0,0)}{i\omega_0}\bigg[\bar{D}\bar{\sigma}e^{-i\omega_0\tau_0}(g_{u_1u_1}(0,0) + i\omega_0g_{u_1u_2}(0,0)) - D\sigma e^{i\omega_0\tau_0}(g_{u_1u_1}(0,0)\\ &- i\omega_0g_{u_1u_2}(0,0))\bigg].	
\end{aligned}
	\end{split}
\end{equation}
Similarly,
\begin{equation} 
	Proj_{{\rm Ker}(M_3^1)}D_xU_2^1(x,0)f_2^1(x,0,0)	= \begin{pmatrix}
		m_{11}x_1^2x_2 + m_{12}x_1x_3^2 \\
		\bar{m}_{11}x_1x_2^2 + \bar{m}_{12}x_2x_3^2 \\
		m_{21}x_1x_2x_3 + m_{22}x_3^3 
		\end{pmatrix},
\end{equation}
where
\begin{equation} 
	\begin{split}
\begin{aligned}
	m_{11} = &\frac{-\tau_0^2\bar{D}\bar{\sigma}}{3i\omega}\bigg[\bar{D}\bar{\sigma}e^{-i\omega_0\tau_0}(-6(g_{u_1u_1}(0,0))^2e^{-i\omega_0\tau_0} + 12i\omega_0g_{u_1u_1}(0,0)g_{u_1u_2}(0,0)e^{-i\omega_0\tau_0}\\ &+ 6\omega_0^4(g_{u_2u_2}(0,0))^2e^{-i\omega_0\tau_0} + 12i\omega_0^3g_{u_1u_2}(0,0)g_{u_2u_2}(0,0)e^{-i\omega_0\tau_0}\\ &-24i\omega_0^3g_{u_1u_2}(0,0)g_{u_2u_2}(0,0) - 24i\omega_0g_{u_1u_1}(0,0)g_{u_1u_2}(0,0)) + D\sigma(14(g_{u_1u_1}(0,0))^2\\ &+ 20\omega_0^2g_{u_1u_1}(0,0)g_{u_2u_2}(0,0) + 14\omega_0^4(g_{u_2u_2}(0,0))^2 + 8\omega_0^2(g_{u_1u_2}(0,0))^2)\\ &+ D_1e^{-i\omega_0\tau_0}(-3(g_{u_1u_1}(0,0))^2 - 3i\omega_0g_{u_1u_2}(0,0)g_{u_2u_2}(0,0) + 3\omega_0^2g_{u_1u_1}(0,0)g_{u_2u_2}(0,0)\\ &- 3i\omega_0g_{u_1u_1}(0,0)g_{u_1u_2}(0,0) - 6\omega_0(g_{u_1u_2}(0,0))^2)\bigg], \\
	m_{12} = &\frac{2\tau_0\bar{D}\bar{\sigma}}{i\omega_0}\bigg[\bar{D}\bar{\sigma}e^{-i\omega_0\tau_0}((g_{u_1u_1}(0,0))^2e^{-i\omega_0\tau_0} + 2i\omega_0g_{u_1u_1}(0,0)g_{u_1u_2}(0,0)\\ &- \omega_0^2g_{u_1u_1}(0,0)g_{u_2u_2}(0,0)e^{-i\omega_0\tau_0})) + D\sigma(-2(g_{u_1u_1}(0,0))^2 - \omega_0^2g_{u_1u_1}(0,0)g_{u_2u_2}(0,0)\\ &- \omega_0^2(g_{u_1u_2}(0,0))^2) + D_1e^{-i\omega_0\tau_0}(2(g_{u_1u_1}(0,0))^2 + 2i\omega_0g_{u_1u_1}(0,0)g_{u_1u_2}(0,0))\bigg], \\
	m_{21} = &\frac{-2\tau_0^2D_1}{i\omega_0}\bigg[\bar{D}\bar{\sigma}e^{-i\omega_0\tau_0}(3(g_{u_1u_1}(0,0))^2 + \omega_0^2g_{u_1u_1}(0,0)g_{u_2u_2}(0,0)\\ &+ 3i\omega_0^3g_{u_1u_2}(0,0)g_{u_2u_2}(0,0) + 2\omega_0^2(g_{u_1u_2}(0,0))^2 + 3i\omega_0g_{u_1u_1}(0,0)g_{u_1u_2}(0,0))\\ &- D\sigma e^{i\omega_0\tau_0}(3(g_{u_1u_1}(0,0))^2 + \omega_0^2g_{u_1u_1}(0,0)g_{u_2u_2}(0,0) - 3i\omega_0^3g_{u_1u_2}(0,0)g_{u_2u_2}(0,0)\\ &+ 2\omega_0^2(g_{u_1u_2}(0,0))^2 - 3i\omega_0g_{u_1u_1}(0,0)g_{u_1u_2}(0,0))\bigg], \\
	m_{22} = &\frac{-2\tau_0^2D_1g_{u_1u_1}(0,0)}{i\omega_0}\bigg[\bar{D}\bar{\sigma}e^{-i\omega_0\tau_0}(g_{u_1u_1}(0,0) + i\omega_0g_{u_1u_2}(0,0)) - D\sigma e^{i\omega_0\tau_0}(g_{u_1u_1}(0,0)\\ &- i\omega_0g_{u_1u_2}(0,0))\bigg].
\end{aligned}
	\end{split}
\end{equation}
And finally, 
\begin{equation}
 	Proj_{{\rm Ker}(M_3^1)}(D_xU_2^1(x,0))^2Jx = \begin{pmatrix}
		n_{11}x_1^2x_2 + n_{12}x_1x_3^2 \\
		\bar{n}_{11}x_1x_2^2 + \bar{n}_{12}x_2x_3^2 \\
		n_{21}x_1x_2x_3 + n_{22}x_3^3 
		\end{pmatrix}, 	
\end{equation}
where
\begin{equation}  
	\begin{split}
\begin{aligned}
	n_{11} = &\frac{-\tau_0^2\bar{D}\bar{\sigma}}{9i\omega_0}\bigg[\bar{D}\bar{\sigma}e^{-i\omega_0\tau_0}(36(g_{u_1u_1}(0,0))^2e^{-i\omega_0\tau_0} + 72i\omega_0g_{u_1u_1}(0,0)g_{u_1u_2}(0,0)\\ &+ 72i\omega_0^3g_{u_1u_2}(0,0)g_{u_2u_2}(0,0) - 36\omega_0^4(g_{u_2u_2}(0,0))^2e^{-i\omega_0\tau_0}) + D\sigma(4(g_{u_1u_1}(0,0))^2\\ &- 8\omega_0^2g_{u_1u_1}(0,0)g_{u_2u_2}(0,0) + 4\omega_0^4(g_{u_2u_2}(0,0))^2 + 16\omega_0^2(g_{u_1u_2})^2)\\ &+D_1e^{-i\omega_0\tau_0}(-9(g_{u_1u_1}(0,0))^2 - 9i\omega_0g_{u_1u_1}(0,0)g_{u_1u_2}(0,0) + 9\omega_0^2g_{u_1u_1}(0,0)g_{u_2u_2}(0,0)\\ &- 9i\omega_0^3g_{u_1u_2}(0,0)g_{u_2u_2}(0,0) - 18\omega_0^2(g_{u_1u_2}(0,0))^2)\bigg] \\
	n_{12} = &\frac{\tau_0^2\bar{D}\bar{\sigma}}{i\omega_0}\bigg[D\sigma((g_{u_1u_1}(0,0))^2 + \omega_0^2(g_{u_1u_2}(0,0))^2) + D_1g_{u_1u_1}(0,0)e^{-i\omega_0\tau_0}(-4g_{u_1u_1}(0,0)\\ &- 4i\omega_0g_{u_1u_2}(0,0))\bigg] \\
	n_{21} = &\frac{\-\tau_0^2D_1}{i\omega_0}\bigg[\bar{D}\bar{\sigma}e^{-i\omega_0\tau_0}((g_{u_1u_1}(0,0))^2 + i\omega_0g_{u_1u_1}(0,0)g_{u_1u_2}(0,0) - \omega_0^2g_{u_1u_1}(0,0)g_{u_2u_2}(0,0)\\ &+ i\omega_0^3g_{u_1u_2}(0,0)g_{u_2u_2}(0,0) + 2\omega_0^2(g_{u_1u_2})^2) - D\sigma e^{i\omega_0\tau_0}(g_{u_1u_1}(0,0))^2\\ &- i\omega_0g_{u_1u_1}(0,0)g_{u_1u_2}(0,0) - \omega_0^2g_{u_1u_1}(0,0)g_{u_2u_2}(0,0) - i\omega_0^3g_{u_1u_2}(0,0)g_{u_2u_2}(0,0)\\ &+ 2\omega_0^2(g_{u_1u_2})^2)\bigg] \\
	n_{22} = &0
\end{aligned}
	\end{split}
\end{equation}

We thus get the desired conclusion.
\hfill\qed

\vspace*{0.25in}
\noindent
{\Large\bf Acknowledgments}

\vspace*{0.2in}
This research is partly supported by the
Natural Sciences and Engineering Research Council of Canada in the
form of a Discovery Grant (BD and VL), and by an Ontario Graduate Scholarship (JB).


\begin{thebibliography}{10}

\bibitem{Atay}
F.M.~Atay.
\newblock Van der Pol's oscillator under delayed feedback.
\newblock {\em J. Sound and Vibration} {\bf 218}, (1998) 333--339.

\bibitem{BBL} 
A.~Beuter, J.~B\'elair and C.~Labrie. 
\newblock Feedback and delays in neurological diseases : a modeling study 
using dynamical systems.
\newblock {\em Bulletin Math. Biology} {\bf 55}, (1993) 525--541.

\bibitem{Cartwright} J. H. E. Cartwright, V. M. Eguiluz, E. Hernandez-Garcia, O. Piro.
\newblock Dynamics of elastic excitable media. 
\newblock {\em Int. J. Bifurcation and Chaos} {\bf 9}, (1999) 2197--2202.

\bibitem{deOlivera}
J.C.F.~de Oliveira.
\newblock Oscillations in a van der Pol equation with delayed argument.
\newblock {\em J. Math. Anal. Appl.} {\bf 275}, (2002) 789--803.

\bibitem{FM1}
T.~Faria and L.T.~Magalh$\tilde{\mbox{\rm a}}$es.
\newblock Normal forms for retarded functional differential equations and applications to Bogdanov-Takens singularity. 
\newblock {\em J. Differential Equations} {\bf 122}, (1995) 201--224.	

\bibitem{FM2}
T.~Faria and L.T.~Magalh$\tilde{\mbox{\rm a}}$es.
\newblock Normal forms for retarded functional differential equations with parameters and applications to Hopf bifurcation.
\newblock {\em J. Differential Equations} {\bf 122}, (1995) 181--200.

\bibitem{GuckHolmes} 
J. Guckenheimer, P. Holmes.
\newblock {\em Nonlinear Oscillations, Dynamical Systems, and Bifurcation of Vector Fields}. 
\newblock Springer-Verlag, New York, 1983.

\bibitem{HVL}
J.K.~Hale and S.M.~Verduyn Lunel.
\newblock Introduction to Functional Differential Equations, 
\newblock Appl. Math. Sci., vol. 99, Springer, New York, 1993.

\bibitem{HFEKGG}
T.~Heil, I.~Fischer, W.~Els\"{a}\mbox{\ss}er, B.~Krauskopf, K.~Green
and A.~Gavrielides.
\newblock Delay dynamics of semiconductor lasers with short external
cavities: Bifurcation scenarios and mechanisms.
\newblock {\em Phys. Rev. E} {\bf 67}, (2003) 066214-1--066214-11.

\bibitem{JiangYuan}
W.~Jiang and Y.~Yuan.
\newblock Bogdanov-Takens singularity in Van der Pol's oscillator with delayed feedback.
\newblock {\em Phys. D} {\bf 227}, (2007) 149--161.

\bibitem{Kaplan} B. Z. Kaplan, I. Gabay, G. Sarafian, D. Sarafian.
\newblock Biological applications of the ÒFilteredÓ Van der Pol oscillator. 
\newblock {\em Journal of the Franklin Institute} {\bf 345}, (2008) 226--232.

\bibitem{Kuang}
Y.~Kuang.
\newblock {\em Delay differential equations with applications in population 
dynamics.} 
\newblock Mathematics in Science and Engineering, 191. Academic Press,  
Boston, (1993).

\bibitem{Kuznetsov} Y. A. Kuznetsov.
\newblock {\em Elements of Applied Bifurcation Theory}.
\newblock Springer-Verlag, New York, 2004, 3rd ed.

\bibitem{LM} 
A.~Longtin and J.G.~Milton. 
\newblock Modeling autonomous oscillations in the human pupil light 
reflex using nonlinear delay-differential equations.
\newblock {\em Bulletin Math. Biology} {\bf 51}, (1989) 605--624.

\bibitem{SieberKrauskopf}
J.~Sieber and B.~Krauskopf.
\newblock Bifurcation analysis of an inverted pendulum with delayed feedback control near a triple-zero eigenvalue singularity.
\newblock {\em Nonlinearity} {\bf 17}, (2004) 85--103.

\bibitem{SC}
E.~Stone and S.A.~Campbell.
\newblock Stability and bifurcation analysis of a nonlinear DDE model
for drilling.
\newblock {\em J. Nonlinear Sci.} {\bf 14}, (2004) 27--57.

\bibitem{SS} M.J. Suarez and P.L. Schopf. 
\newblock A Delayed Action Oscillator for ENSO. 
\newblock
{\em J. Atmos. Sci.}  {\bf 45}, (1988), 3283--3287.

\bibitem{VTK}
A.G.~Vladimirov, D. Turaev and G. Kozyreff.
\newblock Delay differential equations for mode-locked semiconductor
lasers.
\newblock {\em Optics Letters} {\bf 29}, (2004) 1221-1223.

\bibitem{WeiJiang}
J.~Wei and W.~Jiang.
\newblock Stability and bifurcation analysis in Van der Pol's oscillator with delayed feedback.
\newblock {\em J. Sound and Vibration} {\bf 283}, (2005) 801--819.

\bibitem{WeiJiang2} 
J. Wei, W. Jiang.
\newblock
Bifurcation analysis in van der Pol's oscillator with delayed feedback.
\newblock {\em J. Computational and Applied Mathematics} {\bf 213}, (2008) 604--615.

\bibitem{WuWang}
X.~Wu and L.~Wang.
\newblock Zero-Hopf bifurcation for van der Pol's oscillator with delayed feedback.
\newblock {\em J. Comput. Appl. Math.} {\bf 235}, (2011) 2586--2602.

\bibitem{ZW}
C.~Zhang and J.~Wei.
\newblock Stability and bifurcation analysis in a kind of business
cycle model with delay.
\newblock {\em Chaos Solitons Fractals} {\bf 22}, (2004) 883--896.


\end{thebibliography}
\end{document}